\documentstyle[amssymb,10pt]{amsart}
\newtheorem{prop}{Proposition}[section]
\newtheorem{prop:def}{Proposition-Definition}[section]
\newtheorem{defin}{Definition}[section]
\newtheorem{lemma}{Lemma}[section]
\newtheorem{thm}{Theorem}[section]
\newtheorem{cor}{Corollary}[section]
\theoremstyle{remark}

\newtheorem{remark}{Remark}

\begin{document}

\newcommand{\nc}{\newcommand} \nc{\on}{\operatorname}

\nc{\pa}{\partial}

\nc{\cA}{{\cal A}} \nc{\cB}{{\cal B}}\nc{\cC}{{\cal C}} 
\nc{\cD}{{\cal D}} 
\nc{\cE}{{\cal E}} \nc{\cG}{{\cal G}}\nc{\cH}{{\cal H}} 
\nc{\cI}{{\cal I}} \nc{\cJ}{{\cal J}}\nc{\cK}{{\cal K}} 
\nc{\cL}{{\cal L}} 
\nc{\cP}{{\cal P}} \nc{\cQ}{{\cal Q}} 
\nc{\cR}{{\cal R}} \nc{\cS}{{\cal S}}   
\nc{\cV}{{\cal V}} \nc{\cX}{{\cal X}}

\nc{\¦}{{|}}

\nc{\sh}{\on{sh}}\nc{\Id}{\on{Id}}\nc{\Diff}{\on{Diff}}
\nc{\Perm}{\on{Perm}}\nc{\conc}{\on{conc}}\nc{\Alt}{\on{Alt}}
\nc{\ad}{\on{ad}}\nc{\Der}{\on{Der}}\nc{\End}{\on{End}}
\nc{\no}{\on{no\ }} \nc{\res}{\on{res}}\nc{\ddiv}{\on{div}}
\nc{\Sh}{\on{Sh}} \nc{\card}{\on{card}}\nc{\dimm}{\on{dim}}
\nc{\Sym}{\on{Sym}} \nc{\Jac}{\on{Jac}}\nc{\Ker}{\on{Ker}}
\nc{\Spec}{\on{Spec}}\nc{\Cl}{\on{Cl}}
\nc{\Imm}{\on{Im}}\nc{\limm}{\lim}\nc{\Ad}{\on{Ad}}
\nc{\ev}{\on{ev}} \nc{\Hol}{\on{Hol}}\nc{\Det}{\on{Det}}
\nc{\Bun}{\on{Bun}}\nc{\diag}{\on{diag}}\nc{\pr}{\on{pr}} 
\nc{\Span}{\on{Span}}\nc{\Comp}{\on{Comp}}\nc{\Part}{\on{Part}}
\nc{\tensor}{\on{tensor}}\nc{\ind}{\on{ind}}\nc{\id}{\on{id}}
\nc{\Hom}{\on{Hom}}\nc{\Quant}{\on{Quant}} \nc{\Dequant}{\on{Dequant}}\nc{\Def}{\on{Def}}
\nc{\AutLBA}{\on{AutLBA}}\nc{\AutQUE}{\on{AutQUE}}
\nc{\LBA}{{\on{LBA}}}\nc{\Aut}{\on{Aut}}\nc{\QUE}{{\on{QUE}}}
\nc{\QT}{{\on{QT}}}
\nc{\Lyn}{\on{Lyn}}\nc{\Cof}{\on{Cof}}\nc{\LCA}{\underline{\on{LCA}}}\nc{\FLBA}{\on{FLBA}}
\nc{\LA}{\on{LA}}\nc{\FLA}{\on{FLA}}\nc{\EK}{\on{EK}}
\nc{\class}{{\on{class}}}\nc{\br}{\on{br}}\nc{\co}{\on{co}}
\nc{\Prim}{\on{Prim}}\nc{\ren}{\on{ren}}\nc{\lbr}{\on{lbr}}
\nc{\SC}{\on{SC}}\nc{\ol}{\overline}\nc{\FL}{\on{FL}}
\nc{\FA}{\on{FA}}\nc{\alg}{{\on{alg}}}\nc{\KZ}{\on{KZ}}
\nc{\op}{{\on{op}}}\nc{\cop}{{\on{cop}}}
\nc{\Inn}{\on{Inn}}\nc{\OutDer}{\on{OutDer}}
\nc{\inv}{{\on{inv}}}\nc{\gr}{{\on{gr}}}\nc{\Lie}{{\on{Lie}}}
\nc{\Out}{{\on{Out}}}\nc{\univ}{{\on{univ}}}
\nc{\GT}{{\on{GT}}}\nc{\GRT}{{\on{GRT}}}
\nc{\GS}{{\on{GS}}}\nc{\SS}{{\on{SS}}}
\nc{\restr}{{\on{restr}}}\nc{\lin}{{\on{lin}}}
\nc{\mult}{{\on{mult}}}\nc{\qt}{{\on{qt}}}
\nc{\compl}{{\on{compl}}} \nc{\Trees}{{\on{Trees}}}
\nc{\Ord}{{\on{Ord}}} \nc{\norm}{{\on{norm\ ord}}}
\nc{\fd}{{\on{fd}}}\nc{\Maps}{{\on{Maps}}}
\nc{\CYBE}{{\on{CYBE}}} \nc{\CY}{{\on{C}}} \nc{\can}{{\on{can}}}
\nc{\Alg}{{\underline{\on{Alg}}}} \nc{\Poisson}{{\underline{\on{Poisson}}}}
\nc{\Coalg}{{\underline{\on{Coalg}}}} 
\nc{\UE}{{\on{UE}}}\nc{\Bialg}{{\underline{\on{Bialg}}}}
\nc{\QTA}{{\underline{\on{QTA}}}}\nc{\QTQUE}{{\underline{\on{QTQUE}}}}
\nc{\QTLBA}{{\underline{\on{QTLBA}}}}\nc{\TA}{{\underline{\on{TA}}}}
\nc{\MT}{{\underline{\on{MT}}}}\nc{\QYBE}{{\underline{\on{QYBE}}}}
\nc{\QTBialg}{{\underline{\on{QTBialg}}}}
\nc{\qcoco}{{\on{qcoco}}}
\nc{\qcomm}{{\on{qcom}}}\nc{\coco}{{\on{coco}}}
\nc{\dual}{{\on{dual}}} \nc{\sym}{{\on{sym}}}
\nc{\cycl}{{\on{cycl}}} \nc{\Prop}{{\underline{\on{Prop}}}}
\nc{\uV}{{\underline{V}}} \nc{\uR}{{\underline{R}}}
\nc{\graphs}{{\on{graphs}}} \nc{\Isom}{{\on{Isom}}}
\nc{\cp}{{\on{cP}}}\nc{\nat}{{\on{nat}}}\nc{\loc}{{\on{loc}}}
\nc{\PolBMC}{{\on{PolBMC}}}\nc{\coPoisson}{{\on{co-Poisson}}}
\nc{\as}{{\on{antisymm}}}\nc{\Harrison}{{\on{Harrison}}}
\nc{\Mor}{{\on{Mor}}}\nc{\Ob}{{\on{Ob}}}
\nc{\Bil}{{\on{Bil}}} \nc{\invt}{{\on{invt}}}
\nc{\Mod}{{\on{Mod}}} \nc{\quant}{{\on{quant}}}
\nc{\comm}{{\on{comm}}} \nc{\tor}{{\on{tor}}}
\nc{\Irr}{{\on{Irr}}} \nc{\Vect}{{\on{Vect}}}
\nc{\CYBA}{{\on{CYBA}}}\nc{\QYBA}{{\on{QYBA}}}
\nc{\Aug}{{\on{Aug}}} \nc{\Poiss}{{\on{Poiss}}}

\nc{\al}{\alpha}\nc{\g}{\gamma}\nc{\de}{\delta}
\nc{\eps}{\epsilon}\nc{\la}{{\lambda}}
\nc{\si}{\sigma}\nc{\z}{\zeta}

\nc{\La}{\Lambda}

\nc{\ve}{\varepsilon} \nc{\vp}{\varphi} 

\nc{\AAA}{{\mathbb A}}\nc{\BB}{{\mathbb B}}
\nc{\CC}{{\mathbb C}}\nc{\ZZ}{{\mathbb Z}} 
\nc{\QQ}{{\mathbb Q}} \nc{\NN}{{\mathbb N}}\nc{\VV}{{\mathbb V}} 
\nc{\KK}{{\mathbb K}} 

\nc{\ff}{{\mathbf f}}\nc{\bg}{{\mathbf g}}
\nc{\ii}{{\mathbf i}}\nc{\kk}{{\mathbf k}}
\nc{\bl}{{\mathbf l}}\nc{\zz}{{\mathbf z}} 
\nc{\pp}{{\mathbf p}}\nc{\qq}{{\mathbf q}} 

\nc{\cF}{{\cal F}}\nc{\cM}{{\cal M}}\nc{\cO}{{\cal O}}
\nc{\cT}{{\cal T}}\nc{\cW}{{\cal W}}

\nc{\Assoc}{{\mathbf Assoc}}

\def\Sha{{\mathop{\scriptstyle\amalg\!\hspace{-1.8pt}\amalg}}}

\nc{\ub}{{\underline{b}}}
\nc{\uk}{{\underline{k}}} \nc{\ul}{\underline}
\nc{\un}{{\underline{n}}} \nc{\um}{{\underline{m}}}
\nc{\up}{{\underline{p}}}\nc{\uq}{{\underline{q}}}
\nc{\ur}{{\underline{r}}}
\nc{\us}{{\underline{s}}}\nc{\ut}{{\underline{t}}}
\nc{\uw}{{\underline{w}}}
\nc{\uz}{{\underline{z}}}
\nc{\ual}{{\underline{\alpha}}}\nc{\ualpha}{{\underline{\alpha}}}
\nc{\ubeta}{{\underline{\beta}}}\nc{\ugamma}{{\underline{\gamma}}}
\nc{\ueps}{{\underline{\epsilon}}}\nc{\ueta}{{\underline{\eta}}}
\nc{\uzeta}{{\underline{\zeta}}}\nc{\ula}{{\underline{\lambda}}}
\nc{\umu}{{\underline{\mu}}}\nc{\unu}{{\underline{\nu}}}
\nc{\usigma}{{\underline{\sigma}}}\nc{\utau}{{\underline{\tau}}}
\nc{\uI}{{\underline{I}}}\nc{\uJ}{{\underline{J}}}
\nc{\uK}{{\underline{K}}}\nc{\uM}{{\underline{M}}}
\nc{\uN}{{\underline{N}}}

\nc{\A}{{\mathfrak a}} \nc{\B}{{\mathfrak b}} \nc{\C}{{\mathfrak c}} 
\nc{\G}{{\mathfrak g}} \nc{\D}{{\mathfrak d}} \nc{\HH}{{\mathfrak h}}
\nc{\iii}{{\mathfrak i}}\nc{\mm}{{\mathfrak m}}\nc{\N}{{\mathfrak n}} 
\nc{\ttt}{{\mathfrak{t}}}\nc{\U}{{\mathfrak u}}\nc{\V}{{\mathfrak v}}
\nc{\grt}{{\mathfrak grt}}\nc{\gt}{{\mathfrak gt}}
\nc{\SL}{{\mathfrak{sl}}}\nc{\out}{{\mathfrak{out}}}

\nc{\SG}{{\mathfrak S}}

\nc{\wt}{\widetilde} \nc{\wh}{\widehat}

\nc{\bn}{\begin{equation}}\nc{\en}{\end{equation}} \nc{\td}{\tilde}

\title{On the invertibility of quantization functors}

\begin{abstract}
Certain quantization problems are equivalent to 
the construction of morphisms from "quantum" to "classical"
props. Once such a morphism is constructed, Hensel's lemma 
shows that it is in fact an isomorphism. This gives a new, simple 
proof that any Etingof-Kazhdan quantization functor is an  equivalence of 
categories between quantized universal enveloping (QUE) algebras 
and Lie bialgebras over a formal series ring (dequantization).  
We apply the same argument to construct dequantizations of formal 
solutions of the quantum Yang-Baxter equation and of 
quasitriangular QUE algebras. We also give structure results for the
props involved in quantization of Lie bialgebras, which yield an
associator-independent proof that the prop of QUE algebras is a 
flat deformation of the prop of co-Poisson universal enveloping algebras. 
\end{abstract}

\author{Benjamin Enriquez}
\address{IRMA (CNRS et ULP), 7 rue Ren\'e Descartes, F-67084 Strasbourg, France}
\email{enriquez@@math.u-strasbg.fr}

\author{Pavel Etingof}
\address{Department of Mathematics, Massachusetts Institute of Technology,
Cambridge, MA 02139, USA}
\email{etingof@@math.mit.edu}

\maketitle

\section{Introduction}

A prop ("product and permutation category") is an algebraic object 
generalizing the notion of an 
operad (see \cite{McLane}). Given a symmetric monoidal category $\cS$, 
and a prop $\ul P$, one can define the category of $\ul P$-modules 
over $\cS$, $\Mod_{\cS}(\ul P)$. A morphism of props $\ul P\to \ul Q$
then gives rise to a functor $\Mod_{\cS}(\ul Q) \to 
\Mod_\cS(\ul P)$. 

In quantization problems, one should define functors 
from "classical" to "quantum" categories, left inverse to the 
"semiclassical limit" functor. Explicitly, let 
$\cC_\class$ and $\cC_\quant$ be these categories, and 
$\SC : \cC_\quant \to \cC_\class$ be the semiclassical 
limit functor. Then $Q : \cC_\class \to \cC_\quant$ is a 
quantization functor if $\SC \circ Q = \id$. 

In some cases, we have props $\ul P_\class$ and $\ul P_\quant$, 
such that $\cC_x = \Mod_\cS(\ul P_x)$ for $x = \class$ or $\quant$. 
We denote the base field by $\KK$, and by $\hbar$ a formal variable; 
then $\ul P_\quant$ is a module over $\KK[[\hbar]]$, whereas the base ring 
for $\ul P_\class$ is $\KK$. Modules over the prop $\ul P_\quant / (\hbar)$
are provided by $V/(\hbar)$, where $V$ is an object of $\cC_\quant$. 
Such an object carries a classical structure, and is therefore a
$\ul P_\class$-module. This operation has a propic interpretation: 
we have a prop morphism $\ul\SC : \ul P_\class \to \ul P_\quant / (\hbar)$ 
inducing $\SC$. Modules over $\ul P_\class[[\hbar]]$ are provided
by $\hbar$-dependent analogues of the objects of $\cC_\class$; e.g., by the 
$V[[\hbar]]$, where $V\in \Ob(\cC_\class)$ (here the structure maps are 
$\hbar$-independent). 

Then a quantization 
functor $\cC_\class \to \cC_\quant$ may be obtained from a prop morphism 
$\ul Q : \ul P_\quant \to \ul P_\class[[\hbar]]$, such that 
$(\ul Q\ \on{mod}\ \hbar)\circ \ul \SC$
is the identity of $\ul P_\class$. We call such a $\ul Q$ a {\it quantization 
morphism}. 
(Some quantization problems, like quantization of Poisson manifolds or 
algebras, do not fit into this scheme, see Remark \ref{rem:poisson}.)

The main observation of this paper is the following. Assume 
in addition that $\ul\SC$ is surjective. Then Hensel's 
lemma implies that $\ul\SC$ and $\ul Q$ are isomorphisms. Therefore 
the set of quantization morphisms is a torsor, with underlying 
groups $\Aut_1(\ul P_\quant)$ and $\Aut_1(\ul P_\class)$, the subgroups 
of automorphisms of $\ul P_\quant$ and $\ul P_\class[[\hbar]]$ whose reduction 
modulo $\hbar$ is the identity. Moreover, any quantization morphism 
yields an equivalence of categories between $\cC_\quant$ and 
$\cC_{\class,\hbar} = 
\Mod_\cS(\ul P_\class[[\hbar]])$, i.e., between the 
quantum category and the $\hbar$-dependent version 
of the classical category. We call this a dequantization result. 

We apply this to the following three situations: 
(1) quantization of solutions of the CYBE (classical Yang-Baxter 
equation), (2) quantization of Lie bialgebras, (3) quantization 
of quasitriangular Lie bialgebras. Dequantization in situation 
(2) was first obtained in \cite{EK:DQ} using the group $\on{GT}$. 

All three cases are direct applications of the above argument, 
combined in the two last cases with 
the co-Poisson, or quasitriangular versions of the Milnor-Moore theorem. 

In the second situation, we also give an explicit description 
of the structure of the props involved. (A simple description of the 
props involved seems to be impossible in the two other cases.)  
In particular, we prove directly (i.e., not using the existence of 
quantization functors) that the prop $\ul\QUE$ of QUE algebras 
is a flat deformation of the prop $\ul\UE_\cp$ of co-Poisson universal 
enveloping algebras. This implies that any morphism 
$\ul\QUE \to \ul\UE_\cp[[\hbar]]$, whose reduction modulo $\hbar$
is the identity, is an isomorphism. (This by itself does not imply 
the existence of quantization functors, see Remark \ref{rem:poisson}.)

\medskip\noindent
{\bf Acknowledgements.}
The authors thank David Kazhdan, discussions with whom were crucial 
were crucial for many results of this paper. 
P.E. is indebted to IRMA (Strasbourg) for hospitality. 
The research of P.E. was partially supported by
the NSF grant DMS-9988796. 

\section{The formalism of props}

\subsection{Definition, properties}

We fix a base field $\KK$ of characteristic zero, and a 
base ring $R$ containing $\KK$; which will be either $\KK[[\hbar]]$ 
or $\KK$ itself. The modules over $\KK[[\hbar]]$ will always be quotients of 
topologically free modules by closed submodules, and their direct sums 
and tensor products
will be understood in this category; the maps between them will always be 
continuous.  

A prop over $R$ is a symmetric monoidal category 
$\cC$ generated by one object $O$. All the information about such a category 
is contained in the $R$-modules $\Hom_{\cC}(O^{\otimes p},O^{\otimes q}), 
p,q\geq 0$
and the operations relating them.  More specifically, we have: 

\begin{defin} (see \cite{McLane,Lawvere})
A prop $\ul P$ over $R$ is a collection of $R$-modules 
$\ul P(n,m), n,m\geq 0$, together with the data of: 

(1) $R$-module maps $\circ : \ul P(n,m) \otimes \ul P(m,p) \to \ul P(n,p)$
and $\otimes : \ul P(n,m) \otimes \ul P(n',m') \to \ul P(n+n',m+m')$, 
denoted $f\otimes g \mapsto g \circ f$ and 
$f\otimes g \mapsto f \otimes g$,  

(2) linear maps $i_n : \QQ \SG_n \to \ul P(n,n)$, $n\geq 0$, such that 

(a) $\circ$ and $\otimes$ are associative:
$(x\circ y)\circ z = x\circ (y\circ z)$ and 
$(x\otimes y)\otimes z = x\otimes (y\otimes z)$. 
Moreover, we have $(x\circ x')\otimes (y\circ y')
= (x\otimes y) \circ (x'\otimes y')$, 

(b) $i_n$ is an algebra morphism from $\QQ\SG_n$ to 
$(\ul P(n,n),\circ)$, 

(c) for $\sigma\in\SG_n$ and $\sigma'\in\SG_{n'}$, 
denote by $\sigma * \sigma'$ the permutation of $\SG_{n+n'}$
such that $(\sigma * \sigma')(i) = \sigma(i)$ for $i\leq n$ and 
$(\sigma * \sigma')(i) = \sigma'(i - n) + n$ for $i > n$. Then 
$i_{n+n'}(\sigma * \sigma') = i_n(\sigma) \otimes i_{n'}(\sigma')$

(d) if we set $\id = i_1(e)$ ($e$ is the only element of 
$\SG_1$), then the identity $\id^{\otimes m} \circ x = x\circ 
\id^{\otimes n} = x$
holds for $x\in \ul P(n,m)$  

(e) if $\sigma_{n,n'}$ is the permutation in $\SG_{n+n'}$
such that $\sigma_{n,n'}(i) = i+n'$ for $i = 1,\ldots,n$
and  $\sigma_{n,n'}(i) = i-n$ for $i = n+1,\ldots,n+n'$, and
if $x\in \ul P(n,m)$ and $y\in \ul P(n',m')$, then
$$
y\otimes x = \sigma_{m,m'} \circ(x\otimes y) \circ\sigma_{n',n}. 
$$
\end{defin}

If $\cC$ is a symmetric monoidal category generated by $O$, then 
the corresponding prop $\ul P_\cC$ is such that $\ul P_\cC(n,m)
= \Hom_\cC(O^{\otimes n},O^{\otimes m})$. 

If $\ul P$ and $\ul Q$ are two props, then a morphism 
$\phi  : \ul P \to \ul Q$ is a collection of $R$-module
maps $\phi(n,m) : \ul P(n,m) \to \ul Q(n,m)$, such that the 
natural diagrams commute. 

An {\it ideal} $\ul I$ of $\ul P$ is a collection of $R$-submodules 
$\ul I(n,m) \subset \ul P(n,m)$, such that  
$\ul I(n,m) \circ \ul P(m,p) \subset \ul I(n,p)$, 
$\ul P(n,m) \circ \ul I(m,p) \subset \ul I(n,p)$, 
and $\otimes$ takes both $\ul I(n,m) \otimes \ul P(n',m')$ 
and $\ul P(n,m) \otimes \ul I(n',m')$ to $\ul I(n+n',m+m')$. 
The collection of kernels defined by a prop morphism is a prop ideal. 
An ideal $\ul I$ of a prop $\ul P$ gives rise to a quotient 
prop $\ul P / \ul I$, defined by $(\ul P / \ul I)(p,q) = 
\ul P(p,q) / \ul I(p,q)$.

If $\ul P$ be a prop over $R$, then the collection of all 
torsion submodules $\ul P(p,q)_{\on{tor}} \subset \ul P(p,q)$
is an ideal of $\ul P$. We call it the torsion ideal. 

$\ul P$ is a topological prop if it is equipped with a decreasing family 
$\ul I_n$ of prop ideals. We then say that the sequence 
$x_n 
\in \oplus_{p,q}\ul P(p,q)$ tends to zero if $x_n\in 
\oplus_{p,q} \ul I_{k(n)}(p,q)$, where $k(n)$ goes to infinity 
with $n$.  

We will use the following notation. If $x_1,\ldots,x_p$
are such that $x_i\in \ul P(0,n_i)$, if $n = \sum_i n_i$
and $(I_1,\ldots,I_p)$ is a partition of $[1,n]$ by 
ordered sets $I_1,\ldots,I_p$, then $x_1^{I_1} \cdots x_p^{I_p}$
is the element of $\ul P(0,n)$ equal to $\sigma \circ (x_1\otimes \cdots
\otimes x_p)$, where $\sigma\in\SG_n$ is the block permutation  
attached to $I_1,\ldots,I_p$. E.g., $x^{1,4}y^{3,2} = (1432) \circ 
(x\otimes y)$. (We denote by $(i_1\ldots i_k)$ the permutation 
taking $1$ to $i_1$, ..., $k$ to $i_k$.)

\subsection{Props and operads}

Any operad gives rise to a prop. If $(O(n))_{n\geq 0}$
is the family of $\SG_n$-modules underlying an operad, 
then the vector spaces underlying the corresponding prop are 
\begin{equation} \label{rel:prop:operads}
O(n,m) = \bigoplus_{(I_1,\ldots,I_m)\in \Part_m(n)} 
O(\card(I_1)) \otimes\cdots \otimes O(\card(I_m)) . 
\end{equation}
Here $\Part_m(n)$ is the set of partitions of $[1,n]$ by $m$ unordered 
sets. 
So $O(n,m)$ vanishes unless $n\geq m$. A similar construction 
holds with cooperads. 

\subsection{Props defined by generators and relations}

\begin{lemma}
If $V = V(n,m), n,m\geq 0$ is a collection of vector spaces, 
then there is a pair $(\ul P_V,\al_V)$ of a prop $\ul P_V$ and
a collection of linear
maps $\al_{V,n,m} : V(n,m) \to \ul P_V(n,m)$, with the 
following universal property. 
If $(\ul P,\al)$ is any pair of a prop $\ul P$ and 
 a collection of linear maps $\al_{n,m} : 
V(n,m) \to \ul P(n,m)$, then there is a unique prop morphism 
$\al_{\ul P} : \ul P_V\to \ul P$, such that $\al_{\ul P} \circ 
\al_{V} = \al$. 
$\ul P_V$ is unique up to isomorphism, we call it the free prop 
generated by $V$. 
\end{lemma}

{\em Proof.} We construct $\ul P_V$ as follows. Choose a basis 
$(e^{\al}_{i,j})_\al$ of each $V(i,j)$. For each $n,m$, 
let $G_V(n,m)$ be the set of oriented graphs $\Gamma$ of the 
following type. 
Vertices of $\Gamma$  are of three types: "inputs", "outputs"
and "operations". "Operations" vertices correspond to an index 
$(i,j,\al)$. A vertex is said to be of valency $(p,q)$ if it 
has $p$ incoming and $q$ outgoing edges. Input, output 
and $(i,j,\al)$ vertices are of valency $(0,1)$, $(1,0)$
and $(i,j)$. Each vertex carries an order of its input and output  
edges. $\Gamma$ has no oriented cycle. Then $\ul P_V(n,m)$
is the topologically free module spanned by $G_V(n,m)$. 
We define a map $\SG_n \to 
G_V(n,n)$, taking $\sigma$ to the graph of $n$ empty 
edges with (incoming label, outgoing label) 
$= (i,\sigma(i))$. It extends to a linear map $\KK\SG_n\to \ul P_V(n,n)$. 
There are unique maps
$$
\circ^{\on{graphs}} : G_V(n,m) \times G_V(m,p)
\to G_V(n,p)
$$
and 
$$
\otimes^{\on{graphs}} : G_V(n,m) \times G_V(n',m')
\to G_V(n+n',m+m')
$$
defined as follows. If $\Gamma$ and $\Gamma'$ are graphs, 
then $\circ^{\on{graphs}}(\Gamma,\Gamma')$ is obtained from $\Gamma$
and $\Gamma'$ by connecting the output vertex of 
$\Gamma$
with the input vertex of $\Gamma'$ with the same index, and 
then deleting the input and output vertices, and  
$\otimes^{\on{graphs}}(\Gamma,\Gamma')$ is obtained from $\Gamma$
and $\Gamma'$ by adding $n$ (resp., $m$) to the index of each 
input (resp., output) vertex of $\Gamma'$. 
Then $\circ$ and $\otimes$ are 
the linear maps extending  
$\circ^{\graphs}$ and $\otimes^{\graphs}$. 
\hfill \qed\medskip 

Let $V$ be given, and let $\cal R$ be a graded $R$-submodule
of $\oplus_{n,m} \ul P_V(n,m)$. We set  
${\cal R}(n,m) = {\cal R} \cap \ul P_V(n,m)$, so ${\cal R} = 
\oplus_{n,m} {\cal R}(n,m)$. Then we have 

\begin{lemma} 
There exists a unique pair $(\ul P_{V,{\cal R}},\can)$ 
of a prop $\ul P_{V,{\cal R}}$ and a prop morphism 
$\can : \ul P_{V}
\to \ul P_{V,{\cal R}}$, such that $\can({\cal R}) = 0$, 
with the following property. 
If $(\ul Q,\beta)$ is a pair of a prop 
$\ul Q$ and a prop morphism $\beta : \ul P_V \to \ul Q$, such that 
$\beta({\cal R}) = 0$, then there is a unique prop 
morphism $\gamma : \ul P_{V,{\cal R}} \to \ul Q$, such that 
$\gamma\circ\can = \beta$.   
\end{lemma}

{\em Proof.} There is a smallest ideal  
$\ul I_{{\cal R}}$ of $\ul P_V$ (the ideal generated by $\cal R$), 
such that ${\cal R} \subset \oplus_{p,q} \ul I_{\cal R}(p,q)$. 
We then set 
$\ul P_{V,{\cal R}}(n,m) = \ul P_V(n,m) / \ul I_{{\cal R}}(n,m)$. 
\hfill \qed\medskip 

Let us say that two linear combinations of 
graphs of $G_V(n,m)$ are equivalent if 
their difference is a linear combination of 
substitutions of diagrams of $\cal R$ in given graphs. 
Then this equivalence relation is compatible with 
the prop structure, and $\ul P_{V,{\cal R}}(n,m)$ is the quotient of 
$\ul P_V(n,m)$ by this equivalence relation. 

If $\ul P$ is a prop defined by generators and relations, and $\cal R'$
is a collection of new relations involving $x_1,x_2,\ldots$ and the 
generators of 
$\ul P$, we define $\ul P\langle x_1,x_2,\ldots\rangle / (\cal R')$ as the 
prop with generators $\{$generators of $\ul P\}\cup \{x_1,x_2,\ldots\}$
and relations $\{$relations of $\ul P\} \cup \cal R'$ (this definition is 
actually independent on the presentation of $\ul P$). 

\begin{remark}
Any algebra $A$ gives rise to a prop $P_A$, where we  
define $P_A(n,n)$ as the semidirect product of $A^{\otimes n}$
with $\SG_n$, acting on $A^{\otimes n}$ by permutation of factors, and
$P_A(n,m) = 0$ if $n\neq m$; $\circ$ is the product in 
$A^{\otimes n} \rtimes \SG_n$ and $\otimes$ is the product of the 
tensor product and the natural map $\SG_n\times\SG_{n'} \to 
\SG_{n+n'}$. The presentation of a prop by generators and relations
is then a generalization of the similar notion in the case of 
algebras. 
\end{remark}

\subsection{Modules over props}

Let $\cS$ be a symmetric monoidal category over $R$. Then if $A$ is 
any object is $\cS$, the standard operations define a prop 
$\Prop(A)$, where $\Prop(A)(p,q) = \Hom(A^{\otimes p},A^{\otimes q})$. 
A structure of $\ul P$-module over a prop $\ul P$ is a 
pair $(A,\rho)$ of an object $A$ of $\cS$ and a prop morphism 
$\rho : \ul P\to \Prop(A)$. A morphism between two $\ul P$-modules 
$(A, \rho)$ and $(B,\rho')$ is a morphism $\lambda : A\to B$ in $\cS$, 
such that 
if $x\in \ul P(p,q)$ and $a\in A^{\otimes p}$, $\la^{\otimes q}(\rho(x)(a))
= \rho'(x)(\la^{\otimes p}(a))$. 
Then $\ul P$-modules form a category. 

We will sometimes denote $\rho(x)\in \Hom(A^{\otimes p},A^{\otimes q})$
by $x_A$. 

If $\ul P$ is topological, we require that 
the map $\oplus_{p,q} \ul P(p,q) \to  
\End(\oplus_p A^{\otimes p})$ be continuous
in the weak topology: if $x_n\in \oplus_{p,q}\ul P(p,q)$
tends to zero, and $a\in \wh\oplus_{n\geq 0} A^{\otimes n}$, then 
$\rho(x_n)(a)$ tends to zero as $n\to \infty$. Here $\wh\oplus$
denotes the completed direct sum (direct product). When $\cS$ is the 
category of $R$-modules and $R = \KK$, this means 
that $\rho(x_n)(a)$ vanishes for $n$ large enough.

\subsection{Operations on props}

Let $\Irr(n)$ be the set of conjugacy classes of primitive  
idempotents of $\QQ\SG_n$. Let $\pi\in\Irr(n)$, let $\bar\pi\in\QQ\SG_n$
be a representative of $\pi$. The corresponding simple Schur functor 
$F_{(n,\pi)} : \Vect \to \Vect$ is defined by $F_{(n,\pi)}(V) 
= \pi(V^{\otimes n})$. A Schur functor is defined by a multiplicity map 
$\mu : \coprod_{n\geq 0} \Irr(n) \to \ZZ_{\geq 0}$. Then 
$F_\mu(V) := \oplus_{n\geq 0,\pi\in\Irr(n)} F_{(n,\pi)}(V)^{\oplus
\mu(n,\pi)}$. The tensor product of two Schur functors $F,G$ is defined by 
$(F\otimes G)(V) = F(V)\otimes G(V)$. 

If $\ul P$ is a prop and $F = F_\mu$, $F' = F_{\mu'}$ are Schur functors, 
we set 
$$
\ul P(F,F') = \wh\bigoplus_{n\geq 0,\pi\in\Irr(n)}
\bigoplus_{n'\geq 0,\pi'\in\Irr(n')}
(\pi'\circ \ul P(n,n') \circ \pi)^{\oplus \mu(n,\pi)\mu'(n',\pi')}.  
$$

If $F$ is a Schur functor, we define the prop $F(\ul P)$
by $F(\ul P)(p,q) = \ul P(F^{\otimes p},F^{\otimes q})$. 
Then the $F(A)$, $A\in \Mod_{\cS}(\ul P)$, are modules over 
$F(\ul P)$.

\section{Dequantization of solutions of QYBE}
\label{DQ:CYBE}

We denote by $\ul\CYBA$ the prop defined over $\KK$ 
by generators $\eta,m,r$ of bidegrees $(0,1)$, $(2,1)$, $(0,2)$, 
and the following relations: $m \circ (\eta\otimes \id) = m\circ (\id
\otimes \eta) = \id$, $m \circ (m\otimes \id) = m \circ (\id \otimes m)$, 
\begin{equation} \label{propic:CYBE}
(\mu\otimes \id \otimes \id)(r^{1,3}r^{2,4})
+ (\id\otimes\mu \otimes \id)(r^{1,2}r^{3,4}) 
+ (\id\otimes \id\otimes \id)(r^{1,3}r^{2,4}) = 0
\end{equation}
where $\mu = m - m \circ (21)$ 
(the first relations mean that we have a prop morphism $\Alg\to 
\ul\CYBA$, taking $\eta,m$ to their analogues).  

Let $\cS$ be the category of vector spaces, then $\Mod_\cS(\ul\CYBA)$
is the category of quadruples $(A,m_A,1,r_A)$ of an associative algebra 
$(A,m_A,1)$ with unit, together with a solution $r_A\in A^{\otimes 2}$
of the CYBE: $\on{CYB}(r_A) := [r_A^{1,2},r_A^{1,3}] + [r_A^{1,2},r_A^{2,3}]
+ [r_A^{1,3},r_A^{2,3}] = 0$. 

\medskip 

We denote by $\ul\QYBA$ the quotient of the free prop over $\KK[[\hbar]]$
generated by $\eta,m,\rho$ of bidegrees $(0,1)$, $(2,1)$, $(0,2)$, 
by the $\hbar$-adically closed ideal generated by (a) the same
relations as above between $\eta$ and $m$, (b) the relation
\begin{align*}
& (\mu\otimes \id \otimes \id)(\rho^{1,3}\rho^{2,4})
+ (\id\otimes\mu \otimes \id)(\rho^{1,2}\rho^{3,4}) 
+ (\id\otimes \id\otimes \id)(\rho^{1,3}\rho^{2,4}) 
\\ & + \hbar (m\otimes m\otimes m)(\rho^{1,3}\rho^{2,5}\rho^{4,6}
- \rho^{3,5}\rho^{1,6}\rho^{2,4}) 
= 0
\end{align*}
where $\mu = m - m \circ (21)$.  

Let $\cS_\hbar$ be the category of topologically free 
$\KK[[\hbar]]$-modules, 
then $\Mod_{\cS_\hbar}(\ul\QYBA)$ is the category of quadruples 
$(B,m_B,1,\rho_B)$, where $(B,m_B,1)$ is a topologically free algebra, 
together with $\rho_B\in B^{\otimes 2}$, such that 
$\on{CYB}(\rho_B) + \hbar (\rho_B^{1,2}\rho_B^{1,3}\rho_B^{2,3}
- \rho_B^{2,3}\rho_B^{1,3}\rho_B^{1,2})  = 0$. This equation is 
equivalent to the condition that $R_B = 1 + \hbar \rho_B$
satisfies the QYBE (quantum  Yang-Baxter equation). 

\medskip 

Now we have a prop morphism $\ul \SC : \ul\CYBA \to \ul\QYBA/(\hbar)$, 
taking $\eta,m,r$ to the classes of $\eta,m,\rho$. The props $\ul\CYBA$ and 
$\ul\QYBA /(\hbar)$ have the same presentation, therefore $\ul\SC$ is  
an isomorphism. On the other hand, according to \cite{EK,EK:DQ}, 
there exists a prop morphism $\ul Q : \ul\QYBA \to \ul\CYBA[[\hbar]]$, 
such that $(\ul Q\  \on{mod}\ \hbar) \circ \ul\SC$ is the identity. 
This means that $(\ul Q\ \on{mod}\ \hbar)$ is $\ul \SC^{-1}$.  

Recall Hensel's lemma: 

\begin{lemma} If $N$ is a quotient of a topologically free 
$\KK[[\hbar]]$-module by a closed submodule, $M$ is a vector space, 
and $f : N \to M[[\hbar]]$
is a continuous linear map such that $(f\ \on{mod}\ \hbar)$ is 
an isomorphism, 
then $f$ is an isomorphism. In particular, $N$ is torsion-free. 
\end{lemma}

Applying this lemma to the collection of all $\ul\CYBA(p,q)$ and $\ul\QYBA(p,q)$, 
we find: 

\begin{prop} $\ul Q : \ul\QYBA \to \ul\CYBA[[\hbar]]$  is an isomorphism of 
props. 
\end{prop}

Recall that $\ul Q$ takes $m$ to its analogue. For each topologically 
free algebra $(A,m_A,1)$ over $\KK[[\hbar]]$, we have therefore a map 
$r_A\mapsto R(r_A)$
from $\{r_A\in A^{\otimes 2} | r_A$ satisfies the
CYBE$\}$ to $\{R_A\in 1 + \hbar A^{\otimes 2} | R_A$ satisfies the
QYBE$\}$, such that $R(r_A) = 1 + \hbar r_A + O(\hbar^2)$. 
Here $\rho(r_A) =(R(r_A)-1)/\hbar$ is given by a series $r_A +
\sum_{k\geq 2} \hbar^k P_k(m_A,r_A)$, where each $P_k$ a 
certain "polynomial" in 
$m_A,r_A$. For instance, $P_k$ could be equal to 
$\sum_{i,j} a_i a_j b_i \otimes b_j$, 
where $r_A = \sum_i a_i \otimes b_i$.
It is easy to show that such 
a series can be "triangularly" inverted, writing 
$r_A = \rho - \sum_{k\geq 1} \hbar^k P_k(m_A,r_A)$ and substituting 
this expression in this identity iteratively.   

Moreover, we know that the $P_k$ can be chosen to be "normally ordered", 
i.e., in each tensor factor the components $a_i$ are at the left of the 
components $b_j$ (in the language of \cite{Enr:coh}, $P_k$ belongs to 
$(U(\G)^{\otimes 2})_\univ$). 

\begin{cor}
The assignment $r_A\mapsto R(r_A)$
sets up a bijection between $\{$solutions of CYBE in 
$A^{\otimes 2}\}$ and $\{$solutions of QYBE in $1 + 
\hbar A^{\otimes 2}\}$. 
\end{cor}

\section{Dequantization of QUE algebras}

\subsection{The prop $\Bialg$ and related props} 
We denote by $\Bialg$ the prop of bialgebras. It is defined over $\KK$
by generators $m,\Delta,\eta,\epsilon$ of bidegrees 
$(2,1),(1,2),(0,1),(1,0)$ and relations 
$$
m\circ (m\otimes \id) = m\circ (\id\otimes m), \; 
m \circ (\id\otimes \eta) = m \circ (\eta\otimes \id) = \id, 
$$
$$
(\Delta\otimes \id) \circ \Delta = (\id\otimes\Delta) \circ \Delta, \; 
(\eps\otimes\id)\circ \Delta = (\id\otimes\eps)\circ \Delta = \id, 
$$
$$
\Delta \circ m = (m\otimes m) \circ (1324) \circ (\Delta \otimes \Delta), 
\; \Delta \circ \eta = \eta\otimes \eta, \; 
\eps\circ m = \eps\otimes \eps. 
$$
If $\cS$ is the category of $\KK$-vector spaces, then $\Mod_\cS(\Bialg)$
is the category of bialgebras over $\KK$. 

We define the prop $\Bialg_{\qcoco}$ as the quotient of the prop 
$\Bialg[[\hbar]]\langle \wt\delta\rangle / (\Delta - 
(21)\circ \Delta = \hbar \wt\delta)$ by its torsion ideal. 
Here the additional generator $\wt\delta$ has bidegree $(1,2)$. 
If $\cS_\hbar$ is the category of topologically free 
$\KK[[\hbar]]$-modules, then $\Mod_{\cS_\hbar}(\Bialg_{\qcoco})$
is the category of quasi-cocommutative, topologically
free $\KK[[\hbar]]$-bialgebras, i.e., such that $(\Delta - \Delta^{2,1})(A)
\subset \hbar A^{\otimes 2}$. 

We define the prop $\Bialg_{\cp}$ of co-Poisson bialgebras 
as the quotient of $\Bialg\langle \delta\rangle$ ($\delta$ has 
bidegreee $(1,2)$) by the relations
$$
\Delta = (21) \circ \Delta, \; \delta + (21) \circ \delta = 0, \; 
\big( (123) + (231) + (312) \big) \circ (\delta\otimes \id) \circ 
\delta = 0,   
$$
$$
(\Delta \otimes \id)\circ\delta = \big( (123) + (213)\big)
\circ (\id\otimes \delta) \circ \Delta,  
$$  
$$
\delta \circ m = (m\otimes m) \circ (1324) \circ (\delta 
\otimes \Delta + \Delta \otimes \delta), 
$$
$$
\delta \circ \eta = 0, \; (\eps\otimes \id)\circ\delta = 0. 
$$

\begin{lemma}
There is a unique prop morphism $\ul\SC : \Bialg_\cp \to 
\Bialg_{\qcoco} / (\hbar)$, taking $m,\Delta,\eta,\eps$
to their analogues and $\wt\delta$ to $\delta$. $\ul\SC$
is surjective. 
\end{lemma}

{\em Proof.} The proof of the first statement 
is a propic translation of the proof of the 
following fact: if $A$ is a quasi-cocommutative topologically
free $\KK[[\hbar]]$-bialgebra, then 
$A / \hbar A$, equipped with $\delta := 
({{\Delta - \Delta^{2,1}}\over \hbar}\ \on{mod}\ \hbar)$, 
is a co-Poisson bialgebra. Since all generators of 
$\Bialg_{\qcoco} / (\hbar)$ are in the image of $\ul\SC$, 
$\ul\SC$ is surjective. 
\hfill \qed \medskip 

\subsection{Completions} We denote by $\Bialg_\coco$ the quotient of 
$\Bialg$ by the relation $\Delta = (21) \circ \Delta$. Let $I_n$
be the ideal of $\Bialg_{\coco}$ generated by the $(\id - \eta\circ 
\eps)^{\otimes p} \circ \Delta^{(p)}$, $p\geq n$. Here $\Delta^{(n)} = 
(\Delta \otimes \id^{\otimes n-2}) \circ \cdots \circ \Delta$. 
We denote by $\ul\UE$ the completion of $\Bialg_{\coco}$ with respect 
to the family of ideals $I_n$. 

We denote by $J_n$ the ideal of $\Bialg_\cp$ with the same 
generators, and by $\ul\UE_\cp$ the completion of $\Bialg_\cp$
with respect to the family of ideals $J_n$. 

We denote by $K_n$ the ideal of $\Bialg_\qcoco$ with the same generators,
and by $\ul\QUE$ the completion of $\Bialg_\qcoco$ with respect to the family
$K_n$. 

Then $\ul\SC(J_n)$ is contained in the image of $K_n$ under  
$\ul\QUE \to \ul\QUE / (\hbar)$. Therefore: 

\begin{lemma} \label{simple:lemma} \label{lemma:cont}
There is a unique prop morphism 
$\ul\SC : \ul\UE_\cp \to \ul\QUE / (\hbar)$, induced by 
$\ul\SC : \Bialg_\cp \to \Bialg_\qcoco / (\hbar)$, which is 
also surjective. 
\end{lemma}

\subsection{The isomorphism result}

In \cite{EK}, it is shown that there exists a prop morphism 
$\ul Q : \ul\QUE \to \ul\UE_\cp[[\hbar]]$, such that $(\ul Q\ \on{mod}\ \hbar)
\circ \ul\SC = \id$ (i.e., $\ul Q$ is a quantization morphism). 
This implies that $\ul\SC : \ul\UE_\cp \to 
\ul\QUE / (\hbar)$
is injective. Then Lemma \ref{simple:lemma} implies that 
$\ul\SC$ is an isomorphism. 

If now $\ul Q'$ is any quantization morphism, 
then it is a prop morphism 
$\ul Q' : \ul\QUE\to \ul\UE_\cp[[\hbar]]$, such that 
$(\ul Q'\ \on{mod}\ \hbar)$ is an isomorphism. Applying 
Hensel's Lemma to the set of all $\ul Q'(p,q) : \ul\QUE(p,q) \to 
\ul\UE_\cp(p,q)[[\hbar]]$, we get that each $\ul Q'(p,q)$ is an
isomorphism. 

\begin{prop} \label{DQ:morphism}
Each quantization morphism $\ul Q : \ul\QUE \to \ul\UE_\cp[[\hbar]]$
is an isomorphism. So the set of all quantization morphisms
is a torsor over the groups $\Aut_1(\ul\QUE)$ acting on the right 
and $\Aut_1(\ul\UE_\cp[[\hbar]])$ acting on the left. 
\end{prop}


\subsection{Modules over topological props}

Let $\cS$ be the category of vector spaces and let 
us describe the category $\Mod_\cS(\ul\UE)$. 

\begin{lemma}
$\Mod_\cS(\ul\UE)$ is the category of universal enveloping algebras 
over $\KK$, so it is equivalent to the category of Lie algebras. 
\end{lemma}

{\em Proof.} We have a morphism 
$\Bialg_\coco \to \ul\UE$, so if $A$ is a $\ul\UE$-module, then it
is a cocommutative bialgebra $(A,m_A,\Delta_A,\eta_A,\eps_A)$. 
The condition that $A$ is a $\ul\UE$-module means that for
$x_n \in I_n$ and $a\in \oplus_{p} A^{\otimes p}$, 
$\rho(x_n)(a)$ should tend to zero as $n\to\infty$. 
Since the topology of $\oplus_p A^{\otimes p}$ is discrete, 
this means that this sequence vanishes for $n$ large enough. 
In particular, for $a\in A$, the sequence 
$(\id - \eta_A\circ\eps_A)^{\otimes n}
 \circ 
\Delta_A^{(n)}(a)$ vanishes for large $n$. One checks that this 
condition is actually equivalent to $A$ being a $\ul\UE$-module. 
The Milnor-Moore theorem (\cite{MM}) 
then says that $A$ is a universal enveloping 
algebra.\hfill \qed \medskip

It follows that $\Mod_\cS(\ul\UE_\cp)$ is the category 
of universal enveloping algebras with a co-Poisson structure, 
and is equivalent to the category of Lie bialgebras over $\KK$. 

Recall now that $\cS_\hbar$ is the category of topologically
free $\KK[[\hbar]]$-modules. 

\begin{lemma}
$\Mod_{\cS_\hbar}(\ul\UE_\cp[[\hbar]])$ is equivalent 
to the category of topologically free Lie bialgebras over $\KK[[\hbar]]$
(i.e., Lie bialgebras in the category $\cS_\hbar$). 
\end{lemma}

{\em Proof.} The same argument as above shows
that the objects of $\Mod_{\cS_\hbar}(\ul\UE[[\hbar]])$ are the 
topologically free $\KK[[\hbar]]$-bialgebras, such that 
for $a\in A$, the $\hbar$-adic valuation of 
$(\id- \eta_A \circ\eps_A)^{\otimes n} \circ \Delta_A^{(n)}(a)$
tends to zero as $n\to \infty$. A topological version of the 
Milnor-Moore theorem then says that $A$ is the topological 
enveloping algebra of a Lie algebra over $\KK[[\hbar]]$, 
which is a topologically free $\KK[[\hbar]]$-module. 
\hfill \qed \medskip 

We now study $\Mod_{\cS_\hbar}(\ul\QUE)$.

\begin{prop}
The category $\Mod_{\cS_\hbar}(\ul\QUE)$ identifies with the category $\QUE$
of QUE-algebras over $\KK$. 
\end{prop}

{\em Proof.} Let $A$ be a module 
over $\QUE$ in the category $\cS_\hbar$. We have a prop morphism 
$\Bialg_\qcoco\to \ul\QUE$, so $A$ is a quasi-cocommutative 
bialgebra. As above, the condition that $A$ is a $\ul\QUE$-module 
is equivalent to the condition that for each $a\in A$, 
the $\hbar$-adic valuation of 
$(\id - \eta_A \circ \eps_A)^{\otimes n} \circ \Delta_A^{(n)}(a)$
tends to infinity when $n\to\infty$. Let $A_0 = A / \hbar A$. 
Then this condition implies that for each $a_0\in A_0$, 
$(\id - \eta_{A_0} \circ \eps_{A_0})^{\otimes n} \circ \Delta_{A_0}^{(n)}$
vanishes for $n$ large enough. Therefore $A_0$ is a universal 
enveloping algebra. Let us show that $A$ is a Hopf algebra: the antipode 
of $A$ is given by the formula 
$$
S(a) = \eps(a)1 
- a_0 + (a^{(1)})_0 (a^{(2)})_0 - (a^{(1)})_0(a^{(2)})_0(a^{(3)})_0 \ldots, 
$$ 
(we set $x_0 = x - \eps_A(x)1$) i.e., $S = \sum_{n\geq 0} (-1)^{n} m_A^{(n)} \circ 
(\id - \eta_A\circ \eps_A)^{\otimes n} \circ \Delta_A^{(n)}$, 
where $m_A^{(n)}$ is the $n$fold product of $A$. Therefore $A$ is a 
$\KK[[\hbar]]$-Hopf algebra, whose reduction modulo $\hbar$ is a 
universal enveloping algebra, so it is a QUE algebra. 

Conversely, let us show that any QUE algebra $A$ is a $\ul\QUE$-module. 
We should show that for any $a\in A$, the $\hbar$-adic valuation of 
$(\id - \eta_A \circ \eps_A)^{\otimes n} \circ \Delta_A^{(n)}$
tends to infinity with $n$. Let $A_0 = A / \hbar A$. By 
assumption on $A_0$, there exists an integer $n_1$ such that 
$(\id - \eta_{A_0} \circ \eps_{A_0})^{\otimes n_1} \circ \Delta_{A_0}^{(n_1)}
(a\ \on{mod}\ \hbar) = 0$, so 
$(\id - \eta_A \circ \eps_A)^{\otimes n_1} \circ \Delta_A^{(n_1)}
(a) \in \hbar A^{\otimes n_1}$. Let us denote by $a_1$ the class of 
${1\over \hbar}(\id - \eps_A \circ \eta_A)^{\otimes n_1} \circ
\Delta^{(n_1)}(a)$ modulo $\hbar$. This is an element of $A_0^{\otimes n_1}$, 
so there exists $n_2$ such that 
$$
\big( ( (\id - \eta_A \circ \eps_A)^{\otimes n_2} 
\circ \Delta^{(n_2)} )  \otimes \id^{\otimes n_1-1}\big) (a_1) = 0, 
$$
therefore $(\id - \eta_A \circ \eps_A)^{\otimes n_1 + n_2} 
\circ \Delta_A^{(n_1 + n_2)}(a)$ belongs to 
$\hbar^2 A^{\otimes n_1 + n_2}$. In the same way, one constructs a sequence 
of integers $(n_k)_{k\geq 1}$, such that 
$(\id - \eta_A \circ \eps_A)^{\otimes 
n_1 + \cdots + n_k}\circ \Delta_A^{(n_1 + \cdots + n_k)}$ belongs to 
$\hbar^{k} A^{\otimes n_1 + \cdots + n_k}$. This implies that $A$ is a 
$\ul\QUE$-module. \hfill \qed \medskip 

Proposition \ref{DQ:morphism} now implies: 

\begin{thm} (see \cite{EK:DQ})
Each quantization morphism induces an equivalence of categories 
between (a) the category $\QUE$ of QUE-algebras over $\KK$, and 
(b) the category $\on{LBA}_\hbar$ of topologically free 
$\KK[[\hbar]]$-Lie bialgebras. 
\end{thm}

One can define the prop $\ul{\on{Hopf}}$ of Hopf 
algebras as $\ul{\on{Hopf}} = \Bialg\langle S \rangle /$(relations), 
where $S$ has bidegree $(1,1)$ and the relations express the axioms for the 
antipode. Then we have a prop morphism $\ul{\on{Hopf}} \to \ul\QUE$, 
taking $S$ to $\sum_{n\geq 0} (-1)^{n} m^{(n)} \circ 
(\id - \eta\circ\eps)^{\otimes n} \circ \Delta^{(n)}$. We have 
$S^2\in \id + \hbar \ul\QUE(1,1)$, so we have a 1-parameter
subgroup of $\ul\QUE(1,1)^\times$, $\la\mapsto(S^2)^\la$, generated by 
$\on{log}(S^2)\in \ul\QUE(1,1)$. 

\begin{prop} (see \cite{EK3}, Proposition A3) 
Any quantization morphism $\ul\QUE\to\ul\UE_\cp[[\hbar]]$
takes $\on{log}(S^2)$ to a multiple of $\mu\circ\delta$. 
\end{prop}

{\em Proof.} If $\ul P$ is a prop, let us say that a 
prop automorphism of $\ul P$ is $\theta\in \ul P(1,1)$, 
such that $x\circ \theta^{\otimes p} = \theta^{\otimes q}$
for any $x\in \ul P(p,q)$. A prop derivation is the corresponding 
infinitesimal object. Then $S^2$ is a prop automorphism of $\ul\QUE$, 
so if $\ul Q$ is a quantization functor, $\ul Q(S^2)$ is a 
prop automorphism of $\ul\UE_\cp[[\hbar]]$. In particular, is 
commutes with the idempotents $p_n$ (see Lemma \ref{lemma:loday}), 
so it induces a
prop automorphism of $\ul\LBA[[\hbar]]$. Then $\ul Q(\on{\log}(S^2))$
is a prop derivation of $\ul\LBA[[\hbar]]$. In \cite{Enr:univ}, 
we have shown that any such derivation is proportional to 
$\mu\circ\delta$ (here $\mu,\delta$ are the generators of $\ul\LBA$). 
This derivation of $\ul\LBA[[\hbar]]$ extends uniquely 
to a derivation of $\ul\UE_\cp[[\hbar]]$, 
also given by the formula $\mu\circ\delta$  (here $\mu,\delta$ are 
generators of $\ul\UE_\cp$).  
\hfill \qed\medskip 

This proposition was proved in \cite{EK3} when $\ul Q$ is an 
Etingof-Kazhdan quantization morphism. 

\section{Dequantization of QTQUE algebras}

\subsection{Props of some quasitriangular structures}

Recall that the prop $\Bialg_\coco$ is the quotient of the prop 
$\Bialg$ by the ideal generated by $\Delta = (21)\circ \Delta$: 
it is the prop of cocommutative bialgebras. 

Define $\Bialg_{\coco,\qt}$ as $\Bialg_{\coco}\langle r\rangle /$(relations), 
where $r$ has bidegree $(0,2)$ and the relations are: 
$$
(\Delta \otimes \id) \circ r = r^{1,3} \eta^2 + \eta^1 r^{2,3}, \; 
(\id\otimes\Delta) \circ r = r^{1,3} \eta^2 + r^{1,2} \eta^3, \; 
$$
$$
(m\otimes m) \circ (1324) \circ (t\otimes \Delta) = 
(m\otimes m) \circ (1324) \circ (\Delta \otimes t) 
$$
(here $t = r + (21) \circ r$), 
together with the analogue of (\ref{propic:CYBE}). 
The $\Mod_\cS(\Bialg_{\coco,\qt})$ is the category of pairs
$(A,r_A)$, where $A$ is a cocommutative bialgebra, and $r_A\in A^{\otimes 2}$
is such that $(\Delta_A\otimes \id)(r_A) = r_A^{1,3} + r_A^{2,3}$, 
$(\id\otimes\Delta_A)(r_A) = r_A^{1,2} + r_A^{1,3}$,  
the identity $[t_A,\Delta_A(x)] = 0$ holds for any $x\in A$, 
where $t_A = r_A + r_A^{2,1}$, and $\on{CYB}(r_A) = 0$ 
(the two first conditions mean that $r_A\in \on{Prim}(A)^{\otimes 2}$). 
Such a pair $(A,r_A)$ gives rise to a co-Poisson cocommutative bialgebra, 
with $\delta_A(x) = [r_A,\Delta_A(x)]$; this corresponds to a prop morphism 
$\Bialg_\cp \to \Bialg_{\coco,\qt}$. We have an obvious prop 
morphism $\ul\CYBA\to \Bialg_{\coco,\qt}$ (see Section \ref{DQ:CYBE}). 

Define $\Bialg_\QT$ as $\Bialg\langle R,R^{-1}\rangle/$(relations), 
where $R,R^{-1}$ have bidegree $(0,2)$ and the relations are
$$
(\Delta \otimes \id)\circ R = (\id\otimes\id\otimes m)(R^{1,3}R^{2,4}), 
\; 
(\id\otimes\Delta)\circ R = (m\otimes \id\otimes\id)(R^{1,4}R^{2,3}), 
$$
$$
(m\otimes m) \circ (2314) \circ (\Delta\otimes R) = 
(m\otimes m) \circ (1324) \circ (R\otimes \Delta) . 
$$
$$
(m\otimes m)\circ (1324) \circ (R\otimes R^{-1}) = \eta^{\otimes 2}. 
$$

Then  $\Mod_\cS(\Bialg_\QT)$ is the category of quasitriangular 
bialgebras, i.e., pairs $(A,R_A)$, 
where $A$ is a bialgebra and $R_A\in A^{\otimes 2}$ is inverstible, 
such that 
$(\Delta_A\otimes \id)(R_A) = R_A^{1,3}R_A^{2,3}$, 
$(\id\otimes\Delta_A)(R_A) = R_A^{1,3}R_A^{1,2}$, and 
$\Delta_A^{2,1}(x)R_A = R_A \Delta(x)$ holds for any $x\in A$. 
Then $(A,R_A)$ is a solution of the QYBE. 

We define now $\Bialg_{\qcoco,\QT}$ as the quotient of 
$\Bialg_\QT[[\hbar]] \langle \wt\delta,\wt R\rangle / (\Delta 
- (21)\circ\Delta
= \hbar\wt\delta, R = \eta^{\otimes 2} + \hbar \wt R)$ by its torsion ideal.  
Then $\Mod_{\cS_\hbar}(\Bialg_{\qcoco,\QT})$ is the category of 
quasitriangular quasi-cocommutative $\KK[[\hbar]]$-bialgebras
$(A,R_A)$, such that $R_A \in 1 + \hbar A^{\otimes 2}$.  
We have a prop morphism $\ul\QYBA\to \Bialg_{\qcoco,\QT}$, 
taking $m,\eta$ to their analogues and $\rho$ to $\wt R$ 
(see Section \ref{DQ:CYBE}).   

\begin{lemma}
There exists a unique prop morphism $\ul\SC :\Bialg_{\coco,\qt} \to 
\Bialg_{\qcoco,\QT} / (\hbar)$, taking $m,\Delta,\eta,\eps$ to the 
reductions modulo $\hbar$ of their analogues, and taking $r$
to the reduction modulo $\hbar$ of $\wt R$. $\ul\SC$ is surjective.
\end{lemma}

{\em Proof.} The proof that this assignment on generators defines a 
morphism of props is a propic version of the proof of the following fact, 
due to Drinfeld: 
if $(A,R_A)$ is a quasi-cocommutative quasitriangular bialgebra, 
such that $R_A\in 1 + \hbar A^{\otimes 2}$, and if $A_0 = A / \hbar A$, 
$r_A = ({{R_A - 1}\over \hbar}\ \on{mod}\ \hbar)$, then $(A_0,r_A)$
is as above. Let us recall the proof of this fact. The identities 
$(\Delta_A \otimes \id)(R_A) = R_A^{1,3}
R_A^{2,3}$ and $(\id\otimes\Delta_A)(R_A) = R_A^{1,3}R_A^{1,2}$
imply, after we substract $1$, divide by $\hbar$ and reduce modulo 
$\hbar$, that $(\Delta_{A_0}\otimes\id)(r_A) = r_A^{1,3} + 
r_A^{2,3}$ and $(\id\otimes\Delta_{A_0})(r_A) = r_A^{1,2} + r_A^{1,3}$
(in the propic proof, dividing by $\hbar$ is possible because 
$\Bialg_{\qcoco,\QT}$ is constructed to be torsion-free). 
$R_A$ satisfies the QYBE, so substracting from this identity 
$R_A^{1,2} + R_A^{1,3} + R_A^{2,3} -2$, dividing by $\hbar^2$
and reducing modulo $\hbar$, we find that $r_A$ satisfies the CYBE
(again, the propic version uses that $\Bialg_{\qcoco,\QT}$ is 
torsion-free). Finally, we have the identity $(R_A^{2,1}R_A) \Delta_A(x)
= \Delta_A(x) (R_A^{2,1} R_A)$ for any $x\in A$. Substracting 
$\Delta_A(x)$ from both sides, dividing by $\hbar$ and reducing modulo 
$\hbar$, we find that $[r_A + r_A^{2,1},\Delta_{A_0}(x)] = 0$
for any $x\in A_0$ (in the propic case, we use the torsion-freeness 
of $\Bialg_{\qcoco,\QT}$ once more).
All the generators of $\Bialg_{\qcoco,\QT} / (\hbar)$ are in the 
image of $\ul\SC$, so $\ul\SC$ is surjective. 
\hfill \qed \medskip 

\subsection{Completions}

We denote by $I'_n$ the ideal of $\Bialg_{\coco,\qt}$ generated by the 
$(\id - \eta\circ\eps)^{\otimes p} \circ \Delta^{(p)}$, $p\geq n$, 
and $\ul\UE_\qt$ to be the completion of $\Bialg_{\coco,\qt}$ with respect 
to the family $I'_n, n\geq 0$.

We denote by $J'_n$ the ideal of $\Bialg_{\qcoco,\QT}$ generated by the 
analogous elements, and by $\ul\QUE_\QT$ the completion of $\Bialg_{\qcoco,\QT}$
with respect to the family $J'_n,n\geq 0$. Similarly to Lemma \ref{lemma:cont}, 
we have:

\begin{lemma} \label{triv}
$\ul\SC$ extends continuously to a unique morphism of props 
$\ul\SC : \ul\UE_{\coco,\qt}\to \ul\QUE_\QT / (\hbar)$, which is 
also surjective.  
\end{lemma}

\subsection{The isomorphism result}

In \cite{EK:DQ}, it is shown that there exists a prop morphism 
$\ul Q : \ul\QUE_\QT \to \ul\UE_\qt[[\hbar]]$, such that 
$(\ul Q\ \on{mod}\ \hbar) \circ \ul\SC$ is the identity
(i.e., $\ul Q$ is a quantization morphism). 
This implies that $\ul\SC$ is injective. Together with 
Lemma \ref{triv}, this implies that $\ul\SC$ is an 
isomorphism. Now Hensel's lemma implies that any quantization functor 
$\ul Q$ is a prop isomorphism. We have proved: 

\begin{prop} \label{prop:isom:QT}
Any quantization morphism $\ul Q : \ul\QUE_\QT \to \ul\UE_\qt[[\hbar]]$
is a prop isomorphism. 
\end{prop}

\subsection{Modules over quasitriangular props}

\begin{lemma}
$\Mod_\cS(\ul\UE_\qt)$ is equivalent to the category of quasitriangular 
Lie bialgebras over $\KK$, i.e., pairs
$(\A,r_\A)$ of a Lie algebra $\A$ over $\KK$ and $r_\A\in \A^{\otimes 2}$, 
such that $\on{CYB}(r_\A) = 0$ and $r_\A + r_\A^{2,1}$ is invariant. 
\end{lemma}

{\em Proof.} Let $(A,r_A)$ be a module over $\ul\UE_\qt$. 
Then $A$ is a $\ul\UE$-module, so it is a universal enveloping 
algebra. Let $\A = \on{Prim}(A)$. We know that $r_A 
\in \A^{\otimes 2}$, $r_A$ satisfies the CYBE and 
$t_A := r_A + r_A^{2,1}$ commutes with the image of $\Delta : A \to 
A^{\otimes 2}$; therefore $[t_A,x^1 + x^2] = 0$ for any $x\in \A$.
Conversely, since $\A$ generates $A$ as an algebra, 
this last condition implies that $t_A$ commutes with the  
image of $\Delta_A$. \hfill \qed \medskip 

In the same way, one shows that $\Mod_{\cS_\hbar}(\ul\UE_\qt[[\hbar]])$
is equivalent to the category $\LBA_{\qt,\hbar}$ 
of quasitriangular Lie bialgebras  
in the category of topologically free $\KK[[\hbar]]$-modules, 
i.e., of pairs $(\A,r_\A)$, where $\A$ is a topologically free
$\KK[[\hbar]]$-module and $r_\A\in\A^{\otimes 2}$ is a solution of CYBE, 
such that $r_\A + r_\A^{2,1}$ is invariant. 

\begin{lemma}
$\Mod_{\cS_\hbar}(\ul\QUE)$ is equivalent to the category $\QUE_\QT$ of 
quasitriangular QUE algebras over $\KK$, i.e., pairs $(A,R_A)$ of a 
QUE algebra over $\KK$ and $R_A \in 1 + \hbar A^{\otimes 2}$, such that 
$(A,R_A)$ is a 
quasitriangular bialgebra.  
\end{lemma}

These lemmas and Proposition \ref{prop:isom:QT} imply: 

\begin{thm} \label{corr:QT}
Each quantization morphism $\ul Q : \ul\QUE_\QT \to \ul\UE_{\qt}[[\hbar]]$
gives rise to an equivalence of categories between $\QUE_\QT$ and 
$\LBA_{\qt,\hbar}$. 
\end{thm}

Recall that the quantization morphisms from $\cite{EK}$ do not alter the 
algebra structure of $U(\A)$, when $\A$ is quasitriangular. 
When $\ul Q$ is such a quantization 
functor, Theorem \ref{corr:QT} can be made more precise as follows: 

\begin{thm} \label{thm:QT}
Let $\A_0$ be a Lie algebra over $\KK$ and set $\A = \A_0[[\hbar]]$. 
Then each quantization morphism from 
\cite{EK} sets up a bijection between the following coset spaces: 

(a) the set of $r_\A\in \A^{\otimes 2}$, such that $r_\A + r_\A^{2,1}$
is invariant and $\on{CYB}(r_\A) = 0$, modulo the action of 
$\Aut_1(\A)$; 

(b) the set of quasitriangular QUE algebra structures $(\Delta_\A,R_\A)$ 
on $U(\A)$, modulo the action of $\Aut_1(U(\A))$. 

Here $\Aut_1(\A)$ (resp., $\Aut_1(U(\A))$) is the group of 
Lie algebra (resp., algebra) automorphisms
of $\A$ (resp., $U(\A)$), whose reduction modulo $\hbar$ is the identity. 
\end{thm}

{\em Proof.} The proof is based on the following facts: the group 
$\Aut_1(U(\A))$ acts transitively on $\{$cocommutative bialgebra structures on 
$(U(\A),m_0)$ deforming $\Delta_0\}$, by taking $(\theta,\Delta)$ 
to $\theta^{\otimes 2} \circ \Delta \circ \theta^{-1}$. 
The isotropy subgroup of $\Delta_0$ is $\Aut_1(\A)$. 
Here $(m_0,\Delta_0)$ are the undeformed structure maps of $U(\A)$. 
These facts are proved using co-Hochschild 
cohomology.  
\hfill \qed \medskip 

This correspondence is such that $R_\A = 1 + \hbar r_\A + O(\hbar^2)$
and the map $r_\A \mapsto R_\A$ is expressed by the same universal formulas 
as in Section \ref{DQ:CYBE}.

\subsection{}

We will derive from this a classification of twistors related to a given 
associator. 

Let $\A_0$ be a Lie algebra over $\KK$. Let $\A := \A_0[[\hbar]]$
and $t_\A \in S^2(\A)^\A$ be a symmetric invariant tensor. 
Let $\Phi$ be a Drinfeld associator. 
We denote the specialization of $\Phi$ to $(\A,t_\A)$ by 
$\Phi_\A$. 

If $J\in 1 + \hbar U(\A)^{\otimes 2}$, we set 
$$
\wt d(J) := \big(  J^{2,3}J^{1,23} \big)^{-1}
J^{1,2}J^{12,3}. 
$$

In this section, we describe the set $X$ of all 
$J\in 1 + \hbar U(\A)^{\otimes 2}$, such that 
$\wt d(J) = \Phi_\A$. 
We denote by $(u,J) \mapsto u* J$ the action of 
$1 + \hbar U(\A)$ on $1 + \hbar U(\A)^{\otimes 2}$
defined by $u*J := u^1 u^2 J (u^{12})^{-1}$.  
If $\wt d(J)$ is invariant, then $\wt d(u * J) = \wt d(J)$. 

We set $Y = \{\rho\in \A^{\otimes 2} | \on{CYB}(\rho) = 0, 
\rho + \rho^{2,1} = t_\A\}$. 

In \cite{EK,EK:DQ,Enr:coh}, we constructed a map $\rho \mapsto 
J_\Phi(\rho)$, such that if $\rho$ satisfies the CYBE, then 
$\wt d(J_\Phi(\rho)) = \Phi(\hbar \tau^{1,2},\hbar \tau^{2,3})$; 
here $\tau = \rho + \rho^{2,1}$. 
The assignment $\rho\mapsto J_\Phi(\rho)$ defines a map 
from $Y$ to $X$.

Recall that $\exp(\hbar\A)$ is a multiplicative
subgroup of $1 + \hbar U(\A)$. It acts on 
$\{\rho\in \A^{\otimes 2}
| \rho + \rho^{2,1} = t_\A$ and $\on{CYB}(\rho) = 0\}$ 
by conjugation. 

\begin{thm}
Let $J$ be an element of $X$. Then there exists 
$u\in 1 + \hbar U(\A)$ and $\rho\in Y$, such that 
$J = u * J_\Phi(\rho)$. Two pairs $(u,\rho)$ and 
$(u',\rho')$ determine the same $J$ iff there exists an 
element $v\in \exp(\hbar\A)$, such that $u' = uv$ and
$J' = v^{-1} * J$. In other words, $(u,\rho)\mapsto 
u * J_\Phi(\rho)$ defines a bijection 
$$
(1 + \hbar U(\A)) \times_{\exp(\hbar\A)} Y \stackrel{\sim}{\to} X. 
$$
\end{thm}

{\em Proof.} Recall that the Lie algebra 
$\Der(\A,t_\A)$ of derivations of $\A$ leaving
$t_\A$ invariant, acts on $Y$. 

Let $J$ belong to $X$. We will prove the following statement. 
There exist sequences $\rho_n\in Y$, 
$\kappa_n\in U(\A)$, $\gamma_n\in \Der(\A,t_\A)$, 
and algebra automorphisms $\theta_n\in \Aut_1(U(\A))$, such that: 

(1) $J_0 = J$, 
$J_n - J_\Phi(\rho_n) = O(\hbar^n)$, 
$\rho_{n+1} = \exp(\hbar^{n+1}\gamma_{n+1})^{\otimes 2}(\rho_n)$, 
$J_{n+1} = (1 + \hbar^{n+1}\kappa_{n+1}) * J_n$, 
$\theta_n = \id + O(\hbar^n)$; 

(2) $\theta_n$ sets up an isomorphism between the 
quasitriangular QUE algebras 
\begin{equation} \label{QUE1}
(U(\A),m_0, \Ad(J_n) \circ \Delta_0, J_n^{2,1}
e^{\hbar t_\A/2} J_n^{-1}).
\end{equation} 
and 
\begin{equation} \label{QUE2}
(U(\A),m_0,\Ad(J_\Phi(\rho_n)) \circ \Delta_0,
J_\Phi(\rho_n)^{2,1} e^{\hbar t_\A/2} J_\Phi(\rho_n)^{-1}).
\end{equation}
  
We first define $\rho_0$. Twisting the 
quasitriangular quasi-Hopf algebra $(U(\A),m_0,\Delta_0,$ $
e^{\hbar t_\A/2},\Phi_\A)$ by $J$, we obtain a quasitriangular 
QUE algebra 
$$
(U(\A),m_0, \Ad(J) \circ \Delta_0, J^{2,1}
e^{\hbar t_\A/2} J^{-1}).
$$ 
According to Theorem 
\ref{thm:QT}, there exists $\rho_0\in Y$, such that this algebra is 
isomorphic to 
$$
(U(\A),m_0,\Ad(J_\Phi(\rho_0)) \circ \Delta_0,
J_\Phi(\rho_0)^{2,1} e^{\hbar t_\A/2} J_\Phi(\rho_0)^{-1}).
$$ 

We will construct these sequences inductively (the base of induction is 
obvious). 
We will write $J,\rho,J',\rho',\kappa,\theta,\theta''$ instead
of $J_n,\rho_n,J_{n+1}, \rho_{n+1},\kappa_{n+1},\theta_n,\theta_{n+1}$. 

Since the multiplication is the same in algebras (\ref{QUE1})
and (\ref{QUE2}), 
we have $\theta\in\Aut_1(U(\A))$. Moreover, 
\begin{equation} \label{eqs}
\Ad(J) \circ \Delta_0 = \theta^{\otimes 2} \circ \Ad(J_\Phi(\rho))
\circ \Delta_0 \circ \theta^{-1},\; 
J^{2,1} e^{\hbar t_\A/2} J^{-1}
= \theta^{\otimes 2}(
J_\Phi(\rho)^{2,1} e^{\hbar t_\A/2} J_\Phi(\rho)^{-1}). 
\end{equation} 

By hypothesis, we have $J - J_\Phi(\rho)
= O(\hbar^n)$ and $\theta = \id + O(\hbar^n)$. Let 
$K\in U(\A_0)^{\otimes 2}$ and $\gamma \in 
\Der(\A_0,U(\A_0)^{\otimes 2})$ be the reductions modulo 
$\hbar$ of $\hbar^{-n}(J - J_\Phi(\rho))$ and 
$\hbar^{-n}(\theta - \id)$. 
Then (\ref{eqs})  imply 
\begin{equation} \label{eq:K}
[K,\Delta_0(x)] = 
(\gamma\otimes\id + \id \otimes \gamma) (\Delta_0(x)) - \Delta_0(\gamma(x))
\end{equation}
for any $x\in U(\A_0)$, and $K = K^{2,1}$. 
Moreover, since $\wt d(J) = \wt d(J_\Phi(\rho))$, 
we have $\on{d}(K) := K^{12,3} - K^{1,23} - K^{2,3} + K^{1,2} = 0$. 

The equations in $K$ imply that there exists $\kappa\in U(\A_0)$, 
such that $K = \kappa^1 + \kappa^2 - \kappa^{12} =: \on{d}(\kappa)$. 
Then (\ref{eq:K}) implies that for $x\in \A_0$, 
$\gamma(x) - [\kappa,x]\in\A_0$. Therefore $\gamma = \ad(\kappa) + 
\gamma_0$, where $\gamma_0\in \Der(U(\A_0))$ is induced by a 
derivation of $\A_0$. 

We now view $\kappa$ as an element of $U(\A)$
and set $J' := (1 + \hbar^n\kappa)^{-1} * J$. 
Then $J' = J - \hbar^n K + O(\hbar^{n+1})$, therefore 
$J' - J_\Phi(\rho') = O(\hbar^{n+1})$. 
Set $\theta' = \Ad(1 + \hbar^n \kappa)^{-1} \circ \theta$, 
then $\theta' = \id + \hbar^n\gamma_0 + O(\hbar^n)$, where
$\gamma_0$ is viewed as a derivation of $U(\A)$, preserving $\A$. 

Now the second equation in (\ref{eqs}) implies that 
$J' e^{\hbar t_\A} (J')^{-1} $ $= (\theta')^{\otimes 2}
(J_\Phi(\rho) e^{\hbar t_\A} J_\Phi(\rho)^{-1})$. 
So $J' t_\A (J')^{-1} = (\theta')^{\otimes 2}
(J_\Phi(\rho) t_\A J_\Phi(\rho)^{-1})$. The coefficient
of $\hbar^n$ in this identity is yields $(\gamma_0\otimes \id + 
\id\otimes \gamma_0)(t_\A) = 0$. 

Set $\rho' := 
\exp(\hbar^n \gamma_0)^{\otimes 2}(\rho)$. Then 
$J_\Phi(\rho') = J_\Phi(\rho) + O(\hbar^{n+1})$. 
Then $J' - J_\Phi(\rho') = O(\hbar^{n+1})$. 

Set $\theta'' = \theta' \circ 
\exp(\hbar^n \gamma_0)^{-1}$. Then 
$$
\Ad(J') \circ \Delta_0 = \theta^{\prime\prime\otimes 2}
\circ \Ad(J_\Phi(\rho')) \circ \Delta_0 \circ (\theta'')^{-1}, 
\; 
J^{\prime 2,1} e^{\hbar t_\A/2} (J')^{-1} = 
\theta^{\prime\otimes 2}(J_\Phi(\rho')^{2,1}
e^{\hbar t_\A/2} J_\Phi(\rho')^{-1}),  
$$ 
where the second equation follows from the fact that $\gamma_0$
leaves $t_\A$ invariant. These equations mean that 
$\theta'' = \id + O(\hbar^{n+1})$ is an isomorphism between 
the quasitriangular QUE algebras 
$$
(U(\A),m_0,\Ad(J) \circ \Delta_0, J^{2,1} e^{\hbar t_\A/2} J^{-1})
$$
and 
$$
(U(\A),m_0,\Ad(J_\Phi(\rho)) \circ \Delta_0,
J_\Phi(\rho')^{2,1} e^{\hbar t_\A/2} J_\Phi(\rho')^{-1}), 
$$
and we recall that $J' = J_\Phi(\rho') + O(\hbar^{n+1})$. 
This completes the induction step. 

The fact that $u* J_\Phi(\rho) = u' * J_\Phi(\rho')$
implies that $(u,\rho)$ and $(u',\rho')$  are related by the action of 
$\exp(\hbar\A)$ is proved by a co-Hochschild cohomology argument: 
let $n$ be the smallest integer such that $u - u' = O(\hbar^n)$, 
$\rho - \rho' = O(\hbar^n)$. Then if $v,\sigma$ are the reductions modulo 
$\hbar$ of $\hbar^{-n}(u - u')$, $\hbar^{-n}(\rho - \rho')$, then 
$\on{d}(u) + \sigma = 0$, which implies $\sigma = 0$ and 
$u\in \A_0$ by co-Hochschild cohomology.  
\hfill \qed \medskip

\section{Structure results for some props}
\label{sect:str}

\subsection{Props constructed from operads}

We define $\Alg,\Alg_\comm,\Poisson$ as the props associated to the 
operads of associative, commutative and Poisson algebras. 
Let $\cF\cA_N$, $\cF\cC_N$ and $\cF\cP_N$ be the free
associative, commutative and Poisson algebras in $N$ variables
with degrees $\delta_1,\ldots,\delta_N$. Then $\cF\cP_N = S^\cdot 
(\cF\cL_N)$, where $\cF\cL_N$ is the free Lie algebra with $N$
generators, and $\cF\cA_N = U(\cF\cL_N)$.  

Then we have: 

\begin{lemma}
We have 
$$
\Alg(N,n) = (\cF\cA_N^{\otimes n})_{\sum_{i=1}^N \delta_i}, \; 
\Alg_\comm(N,n) = (\cF\cC_N^{\otimes n})_{\sum_{i=1}^N \delta_i}, \; 
\Poisson(N,n) = (\cF\cP_N^{\otimes n})_{\sum_{i=1}^N \delta_i}.  
$$
Here the subscript means the 
homogeneous part of degree $\sum_{i=1}^N \delta_i$. 
\end{lemma}

{\em Proof.} The proof is based on the existence of free 
objects in the categories of associative, Poisson and commutative 
algebras. 
\hfill \qed \medskip 

We can also define the props $\Coalg,\Coalg_\coco,\Coalg_\cp$
of coassociative (resp., cocommutative, co-Poisson) coalgebras,
corresponding to the dual cooperads. Then $\ul{\on{Co}X}(n,N) = \ul X(N,n)$. 

We now define $\Alg_\qcomm$ as the quotient of $\Alg[[\hbar]]\langle 
\wt P\rangle / (m - (21)\circ m = \hbar \wt P)$ by its torsion 
ideal.  Then $\Mod_{\cS_\hbar}(\Alg_\qcomm)$ is the category of 
topologically free, cocommutative $\KK[[\hbar]]$-algebras. 

To describe $\Alg_\qcomm$, we use the following remark. 
Let $M$ be a complete $\KK[[\hbar]]$-module. Let us 
denote by $M_\hbar$ its localization in $\hbar$; this is 
a $\KK((\hbar))$-vector space. Let $M_\tor$ be the torsion 
submodule of $M$, then $M/M_\tor$ is a $\KK[[\hbar]]$-submodule
of $M_\hbar$. 

Let $\Alg((\hbar))$ be the completed tensor product of $\Alg$
with $\KK((\hbar))$, i.e., the version "over $\KK((\hbar))$"
of $\Alg$. Then $\ul U := \Alg((\hbar))\langle \wt P\rangle / 
(m - (21)\circ m = \hbar \wt P)$ coincides with $\Alg((\hbar))$. 
On the other hand, the localization at $\hbar$ of 
$\ul M := \Alg[[\hbar]]\langle \wt P\rangle / 
(m - (21)\circ m = \hbar \wt P)$ coincides with $\ul U$, 
therefore with $\Alg((\hbar))$. So for each $(p,q)$, 
the quotient $\Alg_\qcomm(p,q) = \ul M(p,q) / \ul M(p,q)_\tor$ is a 
$\KK[[\hbar]]$-submodule of $\Alg(p,q)((\hbar))$. 
 
\begin{prop} \label{str:qcomm}
Let $\G$ be the Lie algebra $\cF\cL_N^{\oplus n}$. 
Denote by $U(\G)^{\leq k}$ the linear span of 
the products of less than $k$ elements of $\G$, and by  
$(U(\cF\cL_N)^{\otimes n})^{\leq k}$ the image of $U(\G)^{\leq n}$
under the identification $U(\cF\cL_N)^{\otimes n} = U(\G)$. 
Then we have 
$$
\Alg_{\qcomm}(N,n) \simeq 
\sum_{k\geq 0} \hbar^{k-N} 
\big( U(\cF\cL_N)^{\otimes n} \big)^{\leq k}_{\sum_{i=1}^N \delta_i}
[[\hbar]] 
$$
as a submodule of $\Alg(p,q)((\hbar))$. 
The subscript still denotes the homogeneous 
component of degree $\sum_{i=1}^N \delta_i$. 
\end{prop}

{\em Proof.} Easy. \hfill \qed\medskip 

By construction, $\Alg_{\qcomm}(N,n)$ is a topologically free 
$\KK[[\hbar]]$-module, and $\Alg_\qcomm(N,n) / (\hbar)$
identifies with $\oplus_{k\geq 0} \gr_k(U(\cF\cL_N)^{\otimes n})_{\sum_{i=1}^N 
\delta_i}$, which by the PBW theorem identifies with 
$$
(S^\cdot(\cF\cL_N)^{\otimes n})_{\sum_{i=1}^N \delta_i},
$$ 
i.e., with $\Poisson(N,n)$. 

\begin{cor} \label{cor:poisson}
For any $(N,n)$, we have an isomorphism $\Alg_{\qcomm}(N,n) / (\hbar)
\stackrel{\sim}{\to} \Poisson(N,n)$.  
\end{cor}

We have a morphism $\Poisson\to \Alg_\qcomm / (\hbar)$, taking 
$\eta,m$ to the reduction of their analogues and $P$ to the reduction of 
$\wt P$ (this morphism is a counterpart of 
the functor taking the quasi-commutative algebra $A$ to the 
Poisson algebra $A/\hbar A$). 

Corollary \ref{cor:poisson} then shows that this is an 
isomorphism, so since $\Alg_\qcomm$ is topologically 
free, we get 

\begin{cor} \label{cor:2}
$\Alg_{\qcomm}$ it is a flat deformation of $\Poisson$. 
\end{cor}

\begin{remark} \label{rem:poisson}
Despite this fact, the props
$\Alg_{\qcomm}$ and $\Poisson[[\hbar]]$ are not isomorphic. This 
can be checked explicitly. Besides, it is known that not any Poisson 
algebra can be quantized (see \cite{Mathieu}). 
\end{remark}

\subsection{Props of formal series algebras}

We will say that a formal 
series commutative algebra is an augmented 
commutative algebra $A$, complete for the 
topology defined by the powers of its augmentation ideal $\mm$. 
A formal series Poisson algebra is such an 
algebra, equipped with a Poisson structure $P$, such that 
$P(1,x) = 0$ and $P(x,y)\in \mm$ for any $x,y\in A$.   
Finally, a formal series quasicommutative algebra is a 
quasicommutative augmented algebra over $\KK[[\hbar]]$, 
topologically free as a $\KK[[\hbar]]$-module, complete for the 
topology defined by the powers of $\mm$.  

Define props of augmented commutative (resp., Poisson, 
quasicommutative)  algebras as the props generated by 
$\Alg_\comm$ (resp., $\Poisson$, $\Alg_\qcomm$), the generator
$\eta$ of bidegree $(0,1)$, and the relations 
$\eps\circ\eta = 0$, $\eta\circ m = m \circ (\eps\otimes\eps)$,
together with: in the Poisson case $\eps\circ P = 0$, 
and in the quasicommutative case, $\eps\circ \wt P = 0$. 
We denote them by $\ul\Aug_\comm$, $\ul\Aug_\Poiss$ and $\ul\Aug_\qcomm$. 

Then the corresponding props of formal series algebras are 
defined as the completions of these props with respect to the 
ideals $I^{\comm}_k$, $I_k^\Poiss$ and $I_k^\qcomm$ 
generated by the 
$m^{(l)} \circ (\id - \eta\circ \eps)^{\otimes l}$, $l\geq k$ 
in all three cases. 

Before we describe these ideals, let us describe the props 
$\ul\Aug_X$. 

\begin{lemma}
For any $(N,n)$, the canonical maps followed by composition with 
$(\id - \eta\circ\eps)^{\otimes N}$ induce isomorphisms  
$$
\Alg_\comm(N,n) \simeq
\ul\Aug_\comm(N,n) \circ (\id - \eta\circ\eps)^{\otimes N}, \; 
\Poisson(N,n) \simeq
\ul\Aug_\Poiss(N,n) \circ (\id - \eta\circ\eps)^{\otimes N}, 
$$
$$ 
\Alg_\qcomm(N,n) \simeq
\ul\Aug_\qcomm(N,n) \circ (\id - \eta\circ\eps)^{\otimes N}. 
$$
\end{lemma}
We have therefore identifications 
$$
\ul\Aug_\comm(N,n) \circ (\id - \eta\circ\eps)^{\otimes N}
= (\cF\cC_N^{\otimes n})_{\sum_i\delta_i},  
$$
$$
\ul\Aug_\Poiss(N,n) \circ (\id - \eta\circ\eps)^{\otimes N}
= (S^\cdot(\cF\cL_N)^{\otimes n})_{\sum_i\delta_i},  
$$
$$
\ul\Aug_\qcomm(N,n) 
\circ (\id - \eta\circ\eps)^{\otimes N}
=
\sum_{k\geq 0} \hbar^{k-N} 
(U(\cF\cL_N)^{\otimes n})^{\leq k}_{\sum_i\delta_i}.  
$$

We now describe the intersections of the ideals with these spaces. 
If $\al\geq 1$, we have 
$$
I_\al^\comm(N,n) \circ (\id - \eta\circ\eps)^{\otimes N}
= \oplus_{N\geq \al}(\cF\cC_N^{\otimes n})_{\sum_i\delta_i},  
$$
$$
I_\al^\Poiss(N,n) \circ (\id - \eta\circ\eps)^{\otimes N}
= \oplus_{k\geq \al} (S^k(\cF\cL_N)^{\otimes n})_{\sum_i\delta_i},  
$$
$$
I_\al^\qcomm(N,n)  
\circ (\id - \eta\circ\eps)^{\otimes N}
=
\hbar^{\al - N} (U(\cF\cL_N)^{\otimes n})_{\sum_i\delta_i}[[\hbar]] 
\cap \sum_{k\geq 0} \hbar^{k-N} 
(U(\cF\cL_N)^{\otimes n})^{\leq k}_{\sum_i\delta_i}[[\hbar]].  
$$
One checks that $\ul\Aug_\qcomm$ is a flat deformation of 
$\ul\Aug_\Poiss$ and  $I_\al^\qcomm$ is a flat deformation of $I_\al^P$, 
i.e., it is a saturated subspace whose reduction modulo 
$\hbar$ coincides with $I_\al^\Poiss$.


\subsection{Structures of $\Bialg$ and of the related props}

Let $X$ be one of the indices "no index", cP or coco. 
Then we have prop morphisms $\Alg \to \Bialg_X$ and $\Coalg_X \to 
\Bialg_X$. Composition of these morphisms with the 
operation $\circ$ of $\Bialg_X$ induces linear maps
$\Coalg_X(p,N) \otimes \Alg(N,q) \to \Bialg_X(p,q)$, which 
factor through the natural action of $\SG_N$. 

\begin{prop}
The resulting linear maps
$$
i_{p,q} : \bigoplus_{N\geq 0} \big( 
\Coalg_X(p,N) \otimes \Alg(N,q)\big)_{\SG_N}
\to \Bialg_X(p,q)
$$
are isomorphisms. 
\end{prop} 

{\em Proof.} Let $G$ be a graph for $\Bialg_X(p,q)$. Then the relations 
$\Delta \circ m = (m\otimes m) \circ (1324) \circ (\Delta\otimes \Delta)$, 
together with $\delta\circ m = (m\otimes m) \circ (1324) 
\circ (\delta \otimes \Delta + \Delta \otimes \delta)$
when $X$=cP, 
imply that $G$ can be transformed into a sum of graphs, where
each operation $\Delta$ (and $\delta$ when $X =$ cP) occurs
before each operation $m$. This proves that $i_{p,q}$ is surjective. 

Let us prove that $i_{p,q}$ is injective. The structure of the prop 
$\Alg$ implies that the map 
$$
i : \bigoplus_{N_1,\ldots,N_q\geq 0}
\Coalg_X(p,N_1 + \cdots + N_q)
\to \bigoplus_{N\geq 0} 
\big( 
\Coalg_X(p,N) \otimes \Alg(N,q)\big)_{\SG_N} , 
$$
taking $\oplus_{N_1,\ldots,N_q\geq 0} x_{N_1,\ldots,N_q}$
to $\oplus_{N\geq 0} y_N$, where 
$$
y_N = \sum_{N_1,\ldots,N_q |
\sum_{i=1}^q N_i = N} x_{N_1,\ldots,N_q} \otimes 
(m^{(N_1)} \otimes \cdots \otimes m^{(N_q)}),
$$ 
is a linear isomorphism.
So we should prove that $i_{p,q} \circ i$ is injective. 

Let $\cS = \Vect$. We have a map $\Mod_\cS(\Coalg_X) \to 
\Mod_\cS(\Bialg_X)$, taking a $X$-coalgebra $C$ to $F(C)$. 
Here $F(C)$ is the free associative algebra over the vector space $C$, 
equipped with the unique algebra morphism $\Delta_{F(C)} : 
F(C)\to F(C)^{\otimes 2}$ extending $\Delta: C \to C^{\otimes 2}$, 
and when $X = $cP, with the unique derivation $\delta_{F(C)} : 
F(C) \to F(C)^{\otimes 2}$ extending $\delta_C : C \to C^{\otimes 2}$. 
Let $x = \oplus_{N_1,\ldots,N_q\geq 0} x_{N_1,\ldots,N_q}$, 
then $(i_{p,q}\circ i)(x)_{F(C)}$ is a linear map $F(C)^{\otimes p}
\to F(C)^{\otimes q}$. Composing it to the left with the $p$th 
power of the inclusion $C \hookrightarrow F(C)$ and to the right 
wite the  tensor product of the projections 
$F(C) \to C^{\otimes N_i}$, $i = 1,\ldots,q$, we get a linear map
$C^{\otimes p} \to C^{\otimes N_1 + \cdots + N_q}$, which coincides with 
$(x_{N_1,\ldots,N_q})_C$. This defines a linear map  
$$
\al_C : \Bialg_X(p,q) \to \bigoplus_{N_1,\ldots,N_q\geq 0}
\Hom_\cS(C^{\otimes p},C^{\otimes N_1 + \cdots + N_q}), 
$$
such that $\al_C \circ (i_{p,q} \circ i)$ is the direct sum of the 
prop module maps $\Coalg_X(p,N_1 + \cdots + N_q) \to 
\Hom_\cS(C^{\otimes p},C^{\otimes N_1 + \cdots + N_q})$. 
Taking $C$ to be the cofree $X$-coalgebra with $N_1 + \cdots + N_q$
generators, we see that this map is injective. Therefore $i_{p,q}
\circ i$ is injective. 
\hfill \qed \medskip 

\begin{cor}
Define the linear maps 
$$
j_{p,q} : \bigoplus_{N\geq 0} \big( 
\Coalg_X(p,N) \otimes \Alg(N,q)\big)_{\SG_N}
\to \Bialg_X(p,q)
$$
as the sum the maps taking $x\otimes y$ to 
$\iota_1(y) \circ (\id - \eta\circ\eps)^{\otimes N} \circ \iota_2(x)$, 
where $\iota_1,\iota_2$ are the prop morphisms 
$\Alg\to\Bialg_X$ and $\Coalg_X \to \Bialg_X$. 
Then $j_{p,q}$ is a linear isomorphism. 
\end{cor}

{\em Proof.} Let us denote by $\oplus_{N\geq 0} V_N$ the vector space
on the left. One checks that $j_{p,q} = i_{p,q} \circ k_{p,q}$, 
where $k_{p,q}$ is an endomorphism of $\oplus_{N\geq 0} V_N$, 
whose 
associated graded is the identity for the filtration defined 
by the $\oplus_{N\leq k} V_N$. So $k_{p,q}$ is an isomorphism. 
\hfill \qed \medskip 

The props $\ul\UE$ and $\ul\UE_\cp$ are defined as completions of 
$\Bialg_{\coco}$ and $\Bialg_{\cp}$. We therefore get: 

\begin{prop} \label{str:cp}
The linear maps $j_{p,q}$ extend to linear isomorphisms
$$
\wh j_{p,q} : \wh \bigoplus_{N\geq 0} \big( 
\Coalg_\coco(p,N) \otimes \Alg(N,q)\big)_{\SG_N}
\to \ul\UE(p,q)
$$
and 
$$
\wh j_{p,q} : \wh\oplus_{k\geq 0} \Big( \oplus_{N\geq 0}
\big( \Coalg^k_\cp(p,N)  \otimes \Alg(N,q) \big)_{\SG_N}
\Big) \to \ul\UE_\cp(p,q) , 
$$
where $\Coalg^k_\cp(p,N) = S^k(\cF\cL_N^{\oplus p})_{\sum_i \delta_i}$. 
\end{prop}

The above arguments can be modified to show that 
the analogues of $i_{p,q}$ and $j_{p,q}$ define 
linear isomophisms
$$
\bigoplus_{N\geq 0} \big( 
\Coalg_{\qcoco}(p,N) \otimes \Alg(N,q)\big)_{\SG_N}
\to \Bialg_{\qcoco}(p,q). 
$$ 
We now obtain the structure of the prop $\ul\QUE$. 
We define $\Coalg^{\geq \al}_{\qcoco}(p,N) \subset 
\Coalg_\qcoco(p,N)$ as the 
intersection 
$$
\hbar^{\al - N} (U(\cF\cL_N)^{\otimes p})_{\sum_i\delta_i}[[\hbar]]
\cap 
\sum_{k\geq 0} \hbar^{k-N} 
(U(\cF\cL_N)^{\otimes p})^{\leq k}_{\sum_i\delta_i}[[\hbar]]. 
$$
This space identifies with its dual counterpart $I^\al_\qcomm(N,p)$, 
which is the set of all universally defined 
linear maps $\mm^{\otimes N} \to A^{\otimes p}$ with image 
contained in $(\mm^{(p)})^\al$, where $\mm$ is the augmentation ideal 
of a quasicommutative formal series 
algebra $A$ and $\mm^{(p)}$ is the augmenation 
ideal of $A^{\otimes p}$. 

\begin{prop} $j_{p,q}$ extends to a linear isomorphism 
\begin{align*}
& \wh j_{p,q} : 
\limm_{\leftarrow \al}
\Big( \oplus_{N\geq 0} 
\big( 
\Coalg_{\qcoco}(p,N) \otimes \Alg(N,q)\big)_{\SG_N}
\big/ 
\oplus_{N\geq 0} 
\big( 
\Coalg^{\geq \al}_{\qcoco}(p,N) \otimes \Alg(N,q)\big)_{\SG_N} \big)  
\\ & 
\to \ul\QUE(p,q). 
\end{align*} 
\end{prop}

{\em Proof.} Let $\cI_\al$ be the ideal of 
of $\Bialg_\qcoco$ generated by the $(\id - 
\eta\circ\eps)^{\otimes \beta} \circ\Delta^{(\beta)}$, $\beta\geq \al$. 
Let $\cI_\al^{\on{sat}}$ be its saturation. We have 
$\ul\QUE = \limm_{\leftarrow\al} \Bialg_\qcoco / \cI_\al^{\on{sat}}$. 
 
We should prove that for each $\al$, the ideal $\cI_\al^{\on{sat}}$
is equal to the image 
$\cJ_\al$ of $\oplus_{N\geq 0} 
\big( 
\Coalg^{\geq \al}_{\qcoco}(p,N) \otimes \Alg(N,q)\big)_{\SG_N} \big) $ 
in $\Bialg_\qcoco(p,q)$ under the map $j_{p,q}$. The inclusion 
$\cI_\al^{\on{sat}} \subset \cJ_\al$ is clear, so let 
us show the opposite inclusion. 

Define $\cI'_\al$ as the ideal of $\Bialg$ generated by all 
elements of the form $(\id - \eta\circ\eps)^{\otimes\beta}  \circ \xi$, 
where $\xi\in \Coalg_\qcomm(1,\beta)$ and $\beta\geq \al$. 
Then $\cI'_\al \subset \cI_\al^{\on{sat}}$. We will prove that 
$\cJ_\al \subset \cI'_\al$. 

Set $\delta^{(2)} = (\id - \eta\circ\eps)^{\otimes 2} \circ 
\Delta$ and $\wt \delta^{(2)}
 = (\id - \eta\circ\eps)^{\otimes 2} \circ 
\wt\delta$. Then the key relations are 
\begin{align*}
& \delta^{(2)} \circ m = (m\otimes\id) \circ (132) \circ 
(\delta^{(2)} \otimes\id) +   
(\id\otimes m) \circ 
(\delta^{(2)} \otimes\id) + (m\otimes\id) \circ (\id\otimes\delta^{(2)})
\\ & 
+  (\id\otimes m)\circ (213) \circ (\id\otimes\delta^{(2)}) 
+ (m\otimes m)\circ (1324) \circ (\delta^{(2)}\otimes\delta^{(2)})  
\end{align*}  
and 
\begin{align*}
& \wt \delta^{(2)} \circ m = (m\otimes\id) \circ (132) \circ 
(\wt\delta^{(2)} \otimes\id) +   
(\id\otimes m) \circ 
(\wt\delta^{(2)} \otimes\id) + (m\otimes\id) \circ (\id\otimes\wt\delta^{(2)})
\\ & 
+  (\id\otimes m)\circ (213) \circ (\id\otimes\wt\delta^{(2)}) 
+ (m\otimes m)\circ (1324) \circ (\wt\delta^{(2)}\otimes\delta^{(2)})  
\\ & 
+ (m\otimes m)\circ (1324) \circ \big( \big( (21)\circ \delta^{(2)}\big) 
\otimes
\wt \delta^{(2)}\big).   
\end{align*}  
These relations allow one to show that for any 
$x\in (\id - \eta\circ\eps)^{\otimes \al} \circ 
\Coalg_\qcoco(1,\al)$, $x\circ m$ expressed as a sum 
$\sum_ {\beta,\gamma | \beta\geq \al\on{\ or\ }\gamma\geq \al}
X \otimes (Y\otimes Z)$, where $X\in \Alg(\beta + \gamma,1)$, 
$Y\in (\id - \eta\circ\eps)^{\otimes \beta} \circ 
\Coalg_\qcoco(1,\beta)$
and $Z\in (\id - \eta\circ\eps)^{\otimes \gamma} \circ 
\Coalg_\qcoco(1,\gamma)$. 
These relations allow one to arrange a diagram containing an element of 
$(\id - \eta\circ\eps)^{\otimes \al} \circ 
\Coalg_\qcoco(1,\al)$ as a sum of {\it ordered} diagrams
(i.e., of the form "algebra operations $\circ$ coalgebra
operations"), where all the coalgebra operations are in 
$(\id - \eta\circ\eps)^{\otimes \beta} \circ 
\Coalg_\qcoco(1,\beta)$, $\beta\geq\al$. 
\hfill \qed\medskip

This result, the second part of Proposition \ref{str:cp}, and 
Corollary \ref{cor:2} imply: 

\begin{cor} 
$\ul\QUE$ is a flat deformation of 
$\ul\UE_\cp$. 
\end{cor}

This is a consequence of Proposition \ref{DQ:morphism}, 
but the present proof of this fact 
does not use the existence of quantization functors. 

\begin{remark} Let us describe the spaces $\ul\UE_\cp(1,1)$
and $\ul\QUE(1,1)$. We will view them as spaces of maps 
$\cO\to \cO$, where $\cO$ is a Poisson formal series Hopf algebra
in the first case and a quantized formal series Hopf algebra 
in the second. We define a "Poisson" array of operations 
$$
\begin{array}{cccc}
x & x^{(1)} x^{(2)} & x^{(1)}x^{(2)}x^{(3)} & 
x^{(1)} x^{(2)} x^{(3)} x^{(4)} \\
 & \{x^{(1)} ,x^{(2)} \} & 
\{x^{(1)},x^{(2)}\}x^{(3)} ...
& 
\{x^{(1)} ,x^{(2)} \} x^{(3)}x^{(4)} ... \\
 & & 
 \{\{x^{(1)},x^{(2)}\},x^{(3)}\} ...
& 
\{x^{(1)} ,x^{(2)} \}\{x^{(3)}, x^{(4)}\} ...\\
& & &  \{x^{(1)} ,x^{(2)} \}\{x^{(3)}, x^{(4)} \} ...
\end{array}
$$
and a "quantum" array
$$
\begin{array}{cccc}
f & f^{(1)} f^{(2)} & f^{(1)}f^{(2)}f^{(3)} & 
f^{(1)} f^{(2)} f^{(3)} f^{(4)} \\
 & [f^{(1)} ,f^{(2)} ]_\hbar & 
[f^{(1)},f^{(2)}]_\hbar f^{(3)} ...
& 
[f^{(1)} ,f^{(2)} ]_\hbar f^{(3)}f^{(4)} ... \\
 & & 
 [[f^{(1)},f^{(2)}]_\hbar ,f^{(3)}]_\hbar ...
 & 
\{f^{(1)} ,f^{(2)} ]_\hbar [f^{(3)}, f^{(4)}]_\hbar  ...\\
& & &  [f^{(1)} ,f^{(2)} ]_\hbar [f^{(3)}, f^{(4)}]_\hbar ...
\end{array}
$$
The dots indicate that other monomials belong to a given box, 
and $[-,-]_\hbar = {1\over\hbar}[-,-]$. The bidegree $(i,j)$
"Poisson box" consists of a basis of all polynomials of 
degree $j$ in the $x^{(\al)}$, containing $i-1$ Poisson brackets. 
The Poisson array is graded by the diagonals parallel to 
the main diagonal. The quantum array is filtered by subspaces
lying above the main diagonal (the spaces 
$\Coalg^{\geq \al}_\qcomm(1,1)$). 
An element of $\ul\UE_\cp(1,1)$
is an operation $x\mapsto \sum_{k\geq 0} \sum$
(finite number of elements of the $k$th diagonal), where the 
first sum is infinite.  An element of 
$\ul\QUE(1,1)$ is defined as a similar operation 
$f\mapsto \sum_{k\geq 0} \sum$
(finite number of elements above the $k$th diagonal). 
\hfill \qed\medskip 
\end{remark}

\subsection{The prop of Lie bialgebras and propic Milnor-Moore theorems}

Define $\ul\LBA$ as the prop with generators $\mu,\delta$
with bidegrees $(2,1)$, $(1,2)$, and relations
\begin{equation} \label{rel:mu}
\mu + \mu \circ (21) = 0, \; 
\mu \circ (\mu\otimes \id) \circ ((123) + (231) + (312)) 
 = 0, 
\end{equation}
 \begin{equation} \label{rel:delta} 
\delta + (21) \circ  \mu = 0, \; 
((123) + (231) + (312)) \circ 
(\delta\otimes \id) \circ \delta = 0, 
\end{equation}
$$
\delta \circ \mu = ((12) - (21)) \circ (\id\otimes\mu)\circ 
(\delta\otimes\id) \circ ((12) - (21)). 
$$
Then if $\cS = \Vect$, $\Mod_\cS(\ul\LBA)$ is the category of Lie bialgebras
over $\KK$. 

Define $\ul\LA$ as the prop of Lie algebras, and $\LCA$ as the prop 
of Lie coalgebras. Then $\ul\LA$ is generated by $\mu$ of bidegree 
$(2,1)$ and relations (\ref{rel:mu}), $\LCA$ is generated by 
$\delta$ of bidegree $(1,2)$ and relations (\ref{rel:delta}). 
$\ul\LA$ (resp., $\LCA$) corresponds to the operad (resp., cooperad)
of Lie algebras (resp., coalgebras). We have 
$$
\ul\LA(N,n) = \LCA(n,N) = (\cF\cL_N^{\otimes n})_{\sum_{i=1}^N \delta_i}. 
$$

Moreover, in \cite{Enr:univ,Po}, it is shown that the natural 
prop morphisms $\ul\LA\to\ul\LBA$ and $\LCA \to\ul\LBA$
induce for each $(p,q)$, an isomorphism 
$$
\bigoplus_{N\geq 0} 
(\LCA(p,N) \otimes \ul\LA(N,q))_{\SG_N}
\stackrel{\sim}{\to} \ul\LBA(p,q). 
$$
Therefore, we have an isomorphism 
\begin{equation} \label{form:LBA}
\ul\LBA(p,q) \simeq \bigoplus_{N\geq 0}
\big( (\cF\cL_N^{\otimes p})_{\sum_{i=1}^N\delta_i}
\otimes
(\cF\cL_N^{\otimes q})_{\sum_{i=1}^N\delta_i} \big)_{\SG_N}. 
\end{equation}

On the other hand, we have: 
\begin{thm} \label{propic:MM}
(Propic Milnor-Moore theorem)
We have a prop isomorphism $\ul\UE \stackrel{\sim}{\to} S^\cdot(\ul\LA)$, 
where $S^\cdot = \oplus_{i\geq 0} S^i$ is the "symmetric algebra" Schur
functor. 
\end{thm}

The proof of this theorem is based on the construction of 
"Eulerian idempotents".

\begin{lemma} \label{lemma:proj} \label{lemma:loday} (\cite{Lo})
Let us define the rational numbers $(\la^{(m)}_n)_{n,m\geq 0}$
as the coefficients of the Taylor expansions at zero of 
${1\over{m!}} (\ln(1+u))^m$, so 
$$
{1\over{m!}} (\ln(1+u))^m = \sum_{n\geq 0} \la^{(m)}_n u^n. 
$$
For each $m$, the series
$$
p_m = \sum_{n\geq 0} \la_n^{(m)} m^{(n)} \circ 
(\id - \eta\circ\eps)^{\otimes n}\circ \Delta^{(n)} . 
$$
makes sense in $\ul{\UE}(1,1)$, and the family 
$p_m$ is a complete family of orthogonal idempotents, 
that is $p_m p_{m'} = \delta_{m,m'} p_m$, and 
the sum $\sum_{m\geq 0} p_m$ is equal to $\id$.
\end{lemma}

Moreover, if $\G$ is a Lie algebra, then $U(\G)$ is a $\ul\UE$-module. 
Then $(p_m)_\G\in\End(U(\G))$ corresponds to the projection on the 
$m$th summand of $\oplus_{i\geq 0} S^i(\G)$, under the isomorphism 
$\Sym^{-1} : U(\G) \to S^\cdot(\G)$ (see \cite{Lo}). 

\medskip 

{\em Proof of Theorem \ref{propic:MM}.}
If $p,q,r$ are nonnegative integers, let $\cF\cL_{p+q}$ be the free Lie 
algebra with generators 
$x_1,\ldots,x_p$, $y_1,\ldots,y_q$. Then $\Sym(x_1\otimes\cdots \otimes 
x_p)$ and $\Sym(y_1\otimes\cdots\otimes y_q)$ belong to 
$U(\cF\cL_{p+q})$. So does their product, and it is homogeneous of 
degree 1 in each generator. Let $m_{p,q}^r$ the image of this 
product in $S^r(\cF\cL_{p+q})$ under the composition 
$U(\cF\cL_{p+q}) \stackrel{\Sym^{-1}}{\to} S^\cdot(\cF\cL_{p+q})
\to S^r(\cF\cL_{p+q})$. Then $m_{p,q}^r$ lies in 
$\ul\LA(S^p \otimes S^q,S^r)$, and it vanishes unless $r\leq p+q$. 
Then $m:= \sum_{p,q} \sum_{r=0}^{p+q} m_{p,q}^r$ belongs to 
$S^\cdot(\ul\LA)(1,1)$. We define $\Delta\in S^\cdot(\ul\LA)(1,2)$
by the rule that $\Delta_r^{p,q}$ vanishes unless $r = p+q$, 
and then coincides the propic version of the coproduct for symmetric 
algebras. We define $\eps\in S^\cdot(\ul\LA)(1,0)$ by 
$\eps_i = \delta_{i,0}$ and $\eta\in S^\cdot(\ul\LA)(0,1)$ by 
$\eta^i = \delta_{i,0}$. Then we have a prop morphism 
$\ul\UE \to S^\cdot(\ul\LA)$, taking $m,\delta,\eta,\eps$
to their analogues. 

We now construct a prop morphism $S^\cdot(\ul\LA) \to \ul\UE$. 
Let $p,q$ be integers $\geq 0$, and let $x\in 
S^\cdot(\ul\LA)(p,q)$. We set 
$x = \oplus_{k_i\geq 0,l_i\geq 0} 
x_{k_1,\ldots,k_p}^{l_1,\ldots,l_q}$, where 
$x_{k_1,\ldots,k_p}^{l_1,\ldots,l_q}\in \ul\LA(\otimes_{i=1}^p 
S^{k_i}, \otimes_{j=1}^q S^{l_j})$. We define the map 
$S^\cdot(\ul\LA) \to\ul\UE(p,q)$ to take $x$ to 
\begin{align*}
& \sum_{k_1,\ldots,k_p\geq 0} \sum_{l_1,\ldots,l_q\geq 0}
\big( (m^{(l_1)} \circ \sym^{(l_1)}) 
\otimes \cdots \otimes (m^{(l_q)} \circ \sym^{(l_q)}) \big) 
\circ
\\ & 
\varphi(x_{k_1,\ldots,k_p}^{l_1,\ldots,l_q})
\circ
\big( (p_1^{\otimes k_1} \circ \delta^{(k_1)}) 
\otimes \cdots \otimes (p_1^{\otimes k_p} \circ \delta^{(k_p)}) 
\big) . 
\end{align*}
Here we denote by $\sym_l$ the image of the total symmetrizer
${1\over{l!}}\sum_{\sigma\in\SG_l} \sigma\in \QQ\SG_l$
in $\ul\UE(l,l)$. We denote by $\varphi : \ul\LA\to\ul\UE$
the prop morphism taking $\mu$ to $m - m \circ (21)$. We define 
$\delta^{(p)}\in \ul\UE(1,p)$  as $(\id - \eta\circ\eps)^{\otimes p}
\circ \Delta^{(p)}$. 

This formula corresponds to the following fact. Let $\G$ be a Lie
algebra, and assume that $x$ belongs to $\ul\LA(\otimes_{i=1}^p S^{k_i}, 
\otimes_{j=1}^q S^{l_j})$. Then $x_\G\in \Hom(\otimes_{i=1}^p 
S^{k_i}(\G), \otimes_{j=1}^q S^{l_j}(\G))$. Let $\pi_k : 
U(\G)\to S^k(\G)$ and $i_l : S^l(\G) \to U(\G)$ be the projection 
and injection maps attached to the isomorphism $U(\G) \simeq S^\cdot(\G)$. 
Then 
$$
(\otimes_{j=1}^q i_{l_j}) \circ x_\G \circ (\otimes_{i=1}^p \pi_{k_i}) 
\in \Hom(U(\G)^{\otimes p}, U(\G)^{\otimes q}), 
$$
and it is given by the composition of maps: 
\begin{align*}
& 
U(\G)^{\otimes p} \stackrel{
\otimes_{i=1}^p (p_1^{\otimes k_i} \circ \delta^{(k_i)})_\G}
{-\!\!\!-\!\!\!-\!\!\!-\!\!\!-\!\!\!-\!\!\!-\!\!\!-\!\!\!-\!\!\!-\!\!\!-\!\!\!\!\!
\longrightarrow}
S^{k_1}(\G) \otimes\cdots \otimes S^{k_p}(\G)
\\ & 
\stackrel{x_\G}{\to}
S^{l_1}(\G) \otimes\cdots \otimes S^{l_q}(\G)
\stackrel{
\otimes_{j=1}^q (m^{(l_j)}  \circ \sym_{l_j})_\G 
}
{-\!\!\!-\!\!\!-\!\!\!-\!\!\!-\!\!\!-\!\!\!-\!\!\!-\!\!\!-\!\!\!-
\!\!\!-\!\!\!\!\!\longrightarrow}
U(\G)^{\otimes q}. 
\end{align*}
When writing this diagram, we understand that 
for any $k\geq 0$, $(p_1^{\otimes k} \circ \delta^{(k)})_\G$
maps $U(\G)$ to $S^k(\G)$. The reason why it corresponds to the 
above formula is that if $y\in\ul\LA(p,q)$, then 
$\varphi(y)\in\ul\UE(p,q)$ is such that the restriction of $\varphi(y)_\G$
to $\G^{\otimes p}$ is a map $\G^{\otimes p} \to \G^{\otimes q}$, 
which coincides with $y_\G$. 

One then checks that this is a prop morphism, inverse to 
$\ul\UE \to S^\cdot(\ul\LA)$. 
\hfill \qed \medskip

In the same way, one proves the co-Poisson version of this result: 
\begin{thm}
We have a prop isomorphism $\ul\UE_\cp \stackrel{\sim}{\to}
S^\cdot(\ul\LBA)$, such that the natural diagram 
involving the props $\ul\UE$, $\ul\UE_\cp$, $S^\cdot(\ul\LA)$
and $S^\cdot(\ul\LBA)$ commutes. 
\end{thm}

Taking into account (\ref{form:LBA}), this induces an isomorphism
$$
\ul\UE_\cp(p,q) \simeq \wh\oplus_{k\geq 0}\oplus_{N\geq 0}
\Big( 
\big( S^k(\cF\cL_N^{\oplus p}) \big)_{\sum_{i=1}^N \delta_i}
\otimes 
\big( 
S^\cdot(\cF\cL_N)^{\otimes q}\big)_{\sum_{i=1}^N \delta_i} \Big)_{\SG_N}, 
$$ 
which is the composition with the tensor product of 
$q$ symmetrization maps, of the isomorphism 
$$
\ul\UE_\cp(p,q) 
\simeq \wh\oplus_{k\geq 0} \oplus_{N\geq 0}
\Big( 
\big( S^k(\cF\cL_N^{\oplus p}) \big)_{\sum_{i=1}^N \delta_i}
\otimes 
\big( 
\cF\cA_N^{\otimes q} \big)_{\sum_{i=1}^N \delta_i} \Big)_{\SG_N}
$$ 
given by Proposition \ref{str:cp}.

\end{document}